\documentclass[letterpaper]{article}
\usepackage{theorem,amsmath,amssymb,amscd}
\usepackage[small]{titlesec}

\theoremstyle{change} \allowdisplaybreaks 

\newcommand{\R}{\mathbb{R}}
\newcommand{\C}{\mathbb{C}}
\newcommand{\Z}{\mathbb{Z}}
\newcommand{\Q}{\mathbb{Q}}
\newcommand{\A}{\mathbb{A}}

\newcommand{\SH}{\mathfrak h}
\newcommand{\HH}{\mathbb{H}}
\newcommand{\St}{\mathrm{St}}
\newcommand{\GL}{\mathrm{GL}}
\newcommand{\GU}{\mathrm{GU}}

\newcommand{\SL}{\mathrm{SL}}
\newcommand{\SO}{\mathrm{SO}}
\newcommand{\SU}{\mathrm{SU}}
\newcommand{\GSp}{\mathrm{GSp}}
\newcommand{\SSp}{\mathrm{Sp}}
\newcommand{\OF}{\mathfrak{o}}
\newcommand{\p}{\mathfrak{p}}
\renewcommand{\P}{\mathfrak{P}}

\newcommand{\tr}{{\rm tr}}
\newcommand{\AI}{\mathcal{AI}}
\newcommand{\mat}[4]{{\setlength{\arraycolsep}{0.5mm}\left[
\begin{array}{cc}#1&#2\\#3&#4\end{array}\right]}}
\newcommand{\qed}{\hspace*{\fill}\rule{1ex}{1ex}}
\def\qdots{\mathinner{\mkern1mu\raise0pt\vbox{\kern7pt\hbox{.}}\mkern2mu
\raise3.4pt\hbox{.}\mkern2mu\raise7pt\hbox{.}\mkern1mu}}

\newcommand{\nl}{

\vspace{2ex}}
\newcommand{\nll}{

\vspace{1ex}}

\newtheorem{thm}{Theorem.}[subsection]
\newtheorem{theorem}{Theorem.}[subsection]
\newtheorem{lem}[theorem]{Lemma.}
\newtheorem{lemma}[theorem]{Lemma.}

\newtheorem{proposition}[theorem]{Proposition.}

\begin{document}

\begin{center}
{\Large Integral representation for $L$-functions for $\GSp_4\times \GL_2$, II}

\vspace{2ex} Ameya Pitale\footnote{Department of Mathematics, University of Oklahoma,
Norman, OK 73019-0315, {\tt ameya@math.ou.edu}},
Ralf Schmidt\footnote{Department of Mathematics, University of Oklahoma,
Norman, OK 73019-0315, {\tt rschmidt@math.ou.edu}}

\vspace{3ex}
\begin{minipage}{80ex}\small
{\sc Abstract.} Based on Furusawa's theory \cite{Fu}, we present an integral
representation for the $L$-function $L(s,\pi\times\tau)$, where $\pi$ is a cuspidal
automorphic representation on $\GSp_4$ related to a holomorphic Siegel modular form,
and where $\tau$ is an arbitrary cuspidal automorphic representation on $\GL_2$.
As an application, a special value result for this $L$-function in the spirit
of Deligne's conjecture is proved.
\end{minipage}
\end{center}

\section{Introduction}
Let $F$ be a number field and $\A$ its ring of adeles. Let $\pi$ be a cuspidal,
automorphic representation of $\GSp_4(\A)$, and let $\tau$ be a cuspidal,
automorphic representation of $\GL_2(\A)$. In his paper \cite{Fu}, Furusawa
has obtained an integral representation for the $\GSp_4\times\GL_2$ partial $L$-function
$L^S(s,\pi\times\tau)$ by integrating an Eisenstein series on a unitary group
$\GU(2,2;L)$, where $L$ is a quadratic extension of $F$, against a cusp form
in the space of $\pi$. Furusawa carried out the relevant local $p$-adic
calculations in the case where all the data is unramified, and the relevant
archimedean calculations in the case of ``matching weights''. This was sufficient
to prove a special value result for $L(s,\pi\times\tau)$ in the case where
$\pi$ and $\tau$ come from holomorphic modular forms with respect to the
full modular group.

Furusawa's theory has been refined and extended in the works \cite{PS1} and \cite{Sa}.
In \cite{PS1} it was shown that the $\GL_2$ representation $\tau$ can have
arbitrarily high conductor, as long as its central character remains unramified.
At the archimedean place, the condition on the weights was relaxed. In \cite{Sa} it
was shown that Furusawa's method still works in certain cases of square-free
ramification, both for $\tau$ \emph{and} $\pi$. These more general integral
representations lead to special value results for a wider class of holomorphic
modular forms.

In the present paper we will still assume that $\pi$ is of the type considered
in Furusawa's original work, but we will remove \emph{any} restriction on $\tau$.
Over the number field $\Q$, this means that
$\pi$ is related to a holomorphic Siegel modular form for $\SSp_4(\Z)$, but
we will allow \emph{arbitrary} cuspidal twists $\pi\times\tau$. In a future work
we would like to combine this method with the converse theorem for $\GL_4$
(see \cite{CPS}) in order to lift holomorphic Siegel modular forms (which are non-generic)
to the group $\GL_4$, which is one reason why it is important to be able to twist with
arbitrary cuspidal $\GL_2$ representations.

Via a ``basic identity'' proved in \cite{Fu}, Furusawa's global integrals
factor into an Euler product of local zeta integrals of the form
\begin{equation}\label{intlocalZseq}
 Z(s,W^\#_v,B_v)=\int\limits_{R(F_v)\backslash\GSp_4(F_v)}W_v^\#(\eta h,s)B_v(h)\,dh.
\end{equation}
Here, the function $B_v$ is a vector in a suitable Bessel model of $\pi_v$.
The group $R$ is the corresponding ``Bessel subgroup'' of $\GSp_4$.
The function $W_v^\#$ is a section in a family of induced representations
on the local unitary group $\GU(2,2;L)(F_v)$ (the element $\eta$ is a certain
fixed element in the unitary group). The main point in Furusawa's theory
is to choose the functions $B_v$ and $W^\#_v$ such that the integrals
(\ref{intlocalZseq}) are non-zero for \emph{all} places $v$.

In view of the nature of the representation $\pi$, natural choices of
Bessel functions $B_v$ present themselves, namely as the spherical vector
at finite places, and as a highest weight vector in the archimedean case.
The ``correct'' choice of section $W^\#_v$ is more delicate. For a finite place $v$, the choice of $W^\#$ in \cite{PS1} or \cite{Sa} is a bit {\it ad hoc} and hence not applicable to the case of general $\tau_v$. Here, we will obtain a more natural and canonical choice which works for all representations $\tau_v$ and leads to a much simplified integral calculation in the cases of overlap with \cite{PS1}. We
shall prove that the relevant induced representation on $\GU(2,2;L)(F_v)$ admits
a sort of local newform theory with respect to a certain sequence of
compact-open subgroups (Theorem \ref{unique-W-theorem}). The minimal level
coincides with that of $\tau_v$, the local $\GL_2$ representation from which
the induced representation is constructed. Moreover, at this minimal level,
the space of invariant vectors is one-dimensional. In this sense there is
a unique local newform and we let $W^\#_v$ be this newform. We believe that this choice of local vector in the ramified case is, conceptually, one of the most important contributions of this paper. We hope that it will shed some light on ramified integral calculations in different settings.

Our choice of section in the archimedean case is novel as well,
which allows us to remove a certain assumption on the $\GSp_4$ weight
(condition (4.3.3) in \cite{Fu} and Assumption 2 in \cite{PS1}, 5.1).
But more importantly, the current approach works for non-matching parity
in the $\GSp_4$ and $\GL_2$ weights. Hence, for a real place and suitable
$W^\#_v$, the local integrals (\ref{intlocalZseq}) are always non-zero.
The precise value of the integral is given in Theorem \ref{archmaintheorem},
which is our archimedean main theorem.

In the final part of this paper we
demonstrate how to apply the local theorems in order to derive a global
integral representation for $L(s,\pi\times\tau)$, where the cuspidal automorphic
representation $\pi$ of $\GSp_4(\A_\Q)$ comes from a holomorphic Siegel modular
form, and $\tau$ is an arbitrary cuspidal, automorphic representation of $\GL_2(\A_\Q)$.
Theorem \ref{globalintreptheorem} contains the precise result in the case that
$\tau$ comes from a holomorphic elliptic cusp form of the same weight (even though
such a restriction on $\tau$ is not necessary). We further use this
integral representation to prove a special value result for $L(s,\pi\times\tau)$;
see Theorem \ref{special values thm}. Results of this kind have appeared in
\cite{BH}, \cite{Fu} and \cite{Sa}. The special value result of this paper  substantially adds to the previously known cases, in the sense that it
allows elliptic modular forms (with the same weight) of any level and nebentypus.
In particular, we allow the weights to be odd (the smallest possible odd weight
for a full level cuspidal Siegel modular form is $35$).
 
After some definitions and preliminary remarks, we review Furusawa's general theory
in Sect.\ \ref{furusawa-review}. Following this we develop the non-archimedean theory.
The main result here is the local integral representation
Theorem \ref{ram-central-char-main-thm-local}, but the existence of a ``local newform"
in certain induced representations stated in Theorem \ref{unique-W-theorem} is
possibly of independent interest. For example, the uniqueness of the distinguished
vector is helpful in proving a functional equation, which will be the topic of a
future work. After the non-archimedean theory we develop the archimedean theory,
with the local integral representation Theorem \ref{archmaintheorem} as the
archimedean main result. The final sections contain the global applications
mentioned above.

We would like to thank Abhishek Saha, with whom we had several helpful discussions
on the subject of this paper.
\section{General setup}\label{defn-sec}
In this section, we will recall the basic definitions as stated in Sect.\ 2 of \cite{PS1}.
For simplicity we will make all definitions over a local field, but it is clear
how to define the corresponding global objects.
Let $F$ be a non-archimedean local field of characteristic zero, or
$F=\R$. We fix three elements $a,b,c\in F$ such that
$d:=b^2-4ac\neq0$. Let
\begin{equation}\label{Ldefeq}\renewcommand{\arraystretch}{1.2}
 L = \left\{
      \begin{array}{l@{\qquad\mbox{if }}l}
        F(\sqrt{d})&d \notin F^{\times2},\\
        F \oplus F&d \in F^{\times2}.
      \end{array}
    \right.
\end{equation}
In case $L = F \oplus F$, we consider $F$ diagonally embedded. If
$L$ is a field, we denote by $\bar x$ the Galois conjugate of
$x\in L$ over $F$. If $L=F\oplus F$, let $\overline{(x,y)}=(y,x)$.
In any case we let $N(x)=x\bar x$ and $\tr(x)=x+\bar x$.
\subsection{The unitary group}\label{unitarygroupsec}
We now define the symplectic and unitary similitude groups.
Let $H=\GSp_4$ and $G=\GU(2,2;L)$ be the algebraic $F$-groups whose $F$-points
are given by
\begin{align*}
 H(F)&=\{g \in \GL_4(F) : \, ^tgJg = \mu(g)J,\:\mu(g)\in F^{\times} \}, \\
 G(F)&=\{g \in \GL_4(L) : \, ^t\bar{g}Jg = \mu(g)J,\:\mu(g)\in F^{\times}\},
\end{align*}
where $J = \mat{}{1_2}{-1_2}{}$. Note that $H(F) = G(F)
\cap\GL_4(F)$. As a minimal parabolic subgroup we choose the
subgroup of all matrices that become upper triangular after
switching the last two rows and last two columns. Let $P$ be the
standard maximal parabolic subgroup of $G(F)$ with a non-abelian
unipotent radical. Let $P = MN$ be the Levi decomposition of $P$.
We have $M = M^{(1)}M^{(2)}$, where
\begin{align}
 M^{(1)}(F)&=\{\begin{bmatrix}\zeta\\&1\\&&\bar
  \zeta^{-1}\\&&&1\end{bmatrix}:\:\zeta\in L^\times\} \label{M1defn},\\
 M^{(2)}(F)&=\{\begin{bmatrix}1\\&\alpha&&\beta\\&&\mu\\&\gamma&&\delta\end{bmatrix}:\:
  \alpha,\beta,\gamma,\delta\in L^\times,
  \:\mu=\bar\alpha\delta-\beta\bar\gamma\in F^\times\} \label{M2defn},\\
 N(F) &=\{\begin{bmatrix}
                 1 & z &  &  \\
                  & 1 &  &  \\
                  &  & 1 &  \\
                  &  & -\overline{z} & 1 \\
               \end{bmatrix}
               \begin{bmatrix}
                 1 &  & w & y \\
                  & 1 & \overline{y} &  \\
                  &  & 1 &  \\
                  &  &  & 1 \\
               \end{bmatrix}
              : w\in F,\;y,z \in L \}\label{Ndefn}.
\end{align}
Note that $M^{(2)}(F)\cong\GU(1,1;L)(F)$, where $\GU(1,1;L)$ is defined
analogously to $G=\GU(2,2;L)$. The modular factor of the parabolic $P$ is given by
\begin{equation}\label{deltaPformulaeq}
 \delta_P(\begin{bmatrix}\zeta\\&1\\&&\bar\zeta^{-1}\\&&&1\end{bmatrix}
 \begin{bmatrix}1\\&\alpha&&\beta\\&&\mu\\&\gamma&&\delta\end{bmatrix})
 =|N(\zeta)\mu^{-1}|^3\qquad(\mu=\bar\alpha\delta-\beta\bar\gamma),
\end{equation}
where $|\cdot|$ is the normalized absolute value on $F$.
By Lemma 2.1.1 of \cite{PS1}, the map
\begin{align}\label{M2structurelemmaeq1}
  L^\times\times\GL_2(F)&\longrightarrow\GU(1,1;L)(F),\\
  (\lambda,\mat{\alpha}{\beta}{\gamma}{\delta})&\longmapsto
  \lambda\mat{\alpha}{\beta}{\gamma}{\delta},\nonumber
\end{align}
is surjective with kernel $\{(\lambda,\lambda^{-1}):\:\lambda\in F^\times\}$.
Hence, if $\tau$ is a representation of $\GL_2(F)$, and if $\chi_0$ is a character
of $L^\times$ such that $\chi_0\big|_{F^\times}$ coincides with the central character
of $\tau$, then we can extend $\tau$ to a representation of $\GU(1,1;L)(F)$ by
\begin{equation}\label{M2representationseq}
 \tau(\lambda\mat{\alpha}{\beta}{\gamma}{\delta})
  =\chi_0(\lambda)\tau(\mat{\alpha}{\beta}{\gamma}{\delta}).
\end{equation}
This construction will be used frequently in the following.
Every irreducible, admissible representation of $\GU(1,1;L)(F)$ is of the form
(\ref{M2representationseq}).
\subsection{The Bessel subgroup}\label{besselsubgroupsec}
Recall that we fixed three elements $a,b,c\in F$ such that $d=b^2-4ac\neq0$. Let
$$
 S=\mat{a}{\frac b2}{\frac b2}{c},\qquad
 \xi=\mat{\frac b2}{c}{-a}{\frac{-b}2}.
$$
Then $F(\xi)=F+F\xi$ is a two-dimensional $F$-algebra isomorphic to $L$.
If $L=F(\sqrt{d})$ is a field, then an isomorphism is given by
$x+y\xi\mapsto x+y\frac{\sqrt{d}}2$. If $L=F\oplus F$, then an
isomorphism is given by $x+y\xi\mapsto(x+y\frac{\sqrt{d}}2,x-y\frac{\sqrt{d}}2)$.
The determinant map on $F(\xi)$ corresponds to the norm map on $L$. Let
\begin{equation}\label{TFdefeq}
 T(F)=\{g\in\GL_2(F):\:^tgSg=\det(g)S\}.
\end{equation}
One can check that $T(F)=F(\xi)^\times$. Note that $T(F)\cong
L^\times$ via the isomorphism $F(\xi)\cong L$. We consider $T(F)$
a subgroup of $H(F)$ via
$$
 T(F)\ni g\longmapsto\mat{g}{}{}{\det(g)\,^tg^{-1}}\in H(F).
$$
Let
$$
 U(F)=\{\mat{1_2}{X}{}{1_2}\in H(F):\:^tX=X\}
$$
and $R(F)=T(F)U(F)$. We call $R(F)$ the \emph{Bessel subgroup} of
$H(F)$ (with respect to the given data $a,b,c$). Let $\psi$
be any non-trivial character $F\rightarrow\C^\times$. Let
$\theta:\:U(F)\rightarrow\C^\times$ be the character given by
\begin{equation}\label{thetadefeq}
 \theta(\mat{1}{X}{}{1})=\psi(\tr(SX)).
\end{equation}
Explicitly,
\begin{equation}\label{thetadef2eq}
 \theta(\begin{bmatrix}1&&x&y\\&1&y&z\\&&1\\&&&1\end{bmatrix})=\psi(ax+by+cz).
\end{equation}
We have $\theta(t^{-1}ut)=\theta(u)$ for all $u\in U(F)$ and $t\in T(F)$. Hence,
if $\Lambda$ is any character of $T(F)$, then the map
$tu\mapsto\Lambda(t)\theta(u)$ defines a character of $R(F)$. We denote
this character by $\Lambda\otimes\theta$.
\subsection{Parabolic induction from $P(F)$ to $G(F)$}\label{parabolicinductionsec}
Let $(\tau,V_\tau)$ be an irreducible, admissible representation
of $\GL_2(F)$, and let $\chi_0$ be a character of $L^\times$ such
that $\chi_0\big|_{F^\times}$ coincides with $\omega_{\tau}$, the
central character of $\tau$. Then the pair $(\chi_0,\tau)$ defines a representation
of $M^{(2)}(F)$ as in (\ref{M2representationseq}) on the same space $V_\tau$.
We denote this representation by $\chi_0\times\tau$.
If $V_\tau$ is a space of functions on $\GL_2(F)$ on which
$\GL_2(F)$ acts by right translation, then $\chi_0\times\tau$ can
be realized as a space of functions on $M^{(2)}(F)$ on which
$M^{(2)}(F)$ acts by right translation. This is accomplished by
extending every $W\in V_\tau$ to a function on $M^{(2)}(F)$ via
\begin{equation}\label{extendedWformulaeq}
 W(\lambda g)=\chi_0(\lambda)W(g),\qquad\lambda\in L^\times,\:g\in\GL_2(F).
\end{equation}
If $V_\tau$ is the Whittaker model of $\tau$ with respect to the character $\psi$,
then the extended functions $W$ satisfy the transformation property
\begin{equation}\label{extendedWformulaWhittakereq}
 W(\begin{bmatrix}1\\&1&&x\\&&1\\&&&1\end{bmatrix}g)
 =\psi(x)W(g),\qquad x\in F,\:\:g\in M^{(2)}(F).
\end{equation}
If $s$ is a complex parameter,
$\chi$ is any character of $L^\times$ and $\chi_0\times\tau$ is
a representation of $M^{(2)}(F)$ as above, we denote by $I(s,\chi,\chi_0,\tau)$
the induced representation of $G(F)$ consisting of functions
$f:\:G(F)\rightarrow V_\tau$ with the transformation property
\begin{align}\label{Isfctnspropeq}
 \nonumber &f(\begin{bmatrix}\zeta\\&1\\&&\bar\zeta^{-1}\\&&&1\end{bmatrix}
 \begin{bmatrix}1\\&\lambda\alpha&&\lambda\beta\\
 &&N(\lambda)(\alpha\delta-\beta\gamma)\\&\lambda\gamma&&\lambda\delta\end{bmatrix}ng)\\
 &\hspace{20ex}=\big|N(\zeta\lambda^{-1})(\alpha\delta-\beta\gamma)^{-1}\big|^{3(s+\frac12)}
 \chi(\zeta)\chi_0(\lambda)\tau(\mat{\alpha}{\beta}{\gamma}{\delta})f(g).
\end{align}
Now assume that $V_\tau$ is the Whittaker model of $\tau$ with respect to the
character $\psi$ of $F$. If we associate to each $f$ as above the function
on $G(F)$ given by $W_f(g)=f(g)(1)$, then we obtain another model $I_W(s,\chi,\chi_0,\tau)$
of $I(s,\chi,\chi_0,\tau)$ consisting of functions $W:\:G(F)\rightarrow\C$.
These functions satisfy
\begin{equation}\label{Wsharpproperty1eq}
 W(\begin{bmatrix}\zeta\\&1\\&&\bar\zeta^{-1}\\&&&1\end{bmatrix}
 \begin{bmatrix}1\\&\lambda\\&&N(\lambda)\\&&&\lambda\end{bmatrix}g)
 =|N(\zeta\lambda^{-1})|^{3(s+\frac12)}
 \chi(\zeta)\chi_0(\lambda)W(g),\qquad \zeta,\lambda\in L^\times,
\end{equation}
and
\begin{equation}\label{Wsharpproperty2eq}
 W(\begin{bmatrix}1 & z\\& 1\\&  & 1\\&  & -\overline{z} & 1 \\\end{bmatrix}
  \begin{bmatrix}1&&w&y\\& 1 & \overline{y} &x\\&  & 1\\&  &  & 1 \\\end{bmatrix}g)
 =\psi(x)W(g),\qquad w,x\in F,\;y,z\in L.
\end{equation}
See Sect.\ \ref{furusawa-review} for more on various models for the induced representation.
\subsection{The local integral}\label{localintegralsec}
Let $(\pi,V_\pi)$ be an irreducible, admissible representation of
$H(F)$. Let the Bessel subgroup $R(F)$ be as defined  in
Sect.\ \ref{besselsubgroupsec}; it depends on the given data
$a,b,c\in F$. We assume that $V_\pi$ is a Bessel model for $\pi$
with respect to the character $\Lambda\otimes\theta$ of $R(F)$.
Hence, $V_\pi$ consists of functions $B:\:H(F)\rightarrow\C$
satisfying the Bessel transformation property
$$
 B(tuh)=\Lambda(t)\theta(u)B(h)\qquad
 \text{for }t\in T(F),\:u\in U(F),\:h\in H(F).
$$
Let $(\tau,V_\tau)$ be a generic, irreducible, admissible
representation of $\GL_2(F)$ such that $V_\tau$ is the
$\psi^{-c}$--Whittaker model of $\tau$ (we assume $c\neq0$). Let
$\chi_0$ be a character of $L^\times$ such that
$\chi_0\big|_{F^\times}=\omega_\tau$. Let $\chi$ be the character
of $L^\times$ which satisfies
\begin{equation}\label{chi-lambda-char-condition}
 \chi(\zeta) = \Lambda(\bar{\zeta})^{-1} \chi_0(\bar{\zeta})^{-1}.
\end{equation}
Let $W^\#(\,\cdot\,,s)$ be an element of $I_W(s,\chi,\chi_0,\tau)$ for
which the restriction of $W^\#(\,\cdot\,,s)$ to the standard
maximal compact subgroup of $G(F)$  is
independent of $s$, i.e., $W^\#(\,\cdot\,,s)$ is a ``flat section'' of
the family of induced representations $I_W(s,\chi,\chi_0,\tau)$. By
Lemma 2.3.1 of \cite{PS1}, it is meaningful to consider the
integral
\begin{equation}\label{localZseq}
 Z(s,W^\#,B)=\int\limits_{R(F)\backslash H(F)}W^\#(\eta h,s)B(h)\,dh.
\end{equation}
Here,
\begin{equation}\label{alphadefeq}
 \eta=\begin{bmatrix}1\\\alpha&1\\&&1&-\bar\alpha\\&&&1\end{bmatrix},
 \qquad\text{where}\qquad
 \alpha=\left\{\begin{array}{l@{\qquad\text{if }L}l}
 \displaystyle\frac{b+\sqrt{d}}{2c}&\text{ is a field},\\[2ex]
 \displaystyle\Big(\frac{b+\sqrt{d}}{2c},\frac{b-\sqrt{d}}{2c}\Big)&=F\oplus F.
 \end{array}\right.
\end{equation}
As explained in Sect.\ \ref{furusawa-review} below, see in particular (\ref{localzetaeq}),
these local integrals appear in integral representations
for global $\GSp_4\times\GL_2$ $L$-functions. Therefore,
being able to make choices for the functions $W^\#$ and $B$ such that $Z(s,W^\#,B)$ is
non-zero leads to such integral representations.
In the following we shall demonstrate that this is always possible for local
$\GSp_4(F)$ representations $\pi$ that are relevant for the global application
to Siegel modular forms we have in mind. In the real case we shall
assume that $\pi$ is a holomorphic discrete series representation
and that $B$ corresponds to the highest weight vector. In the
$p$-adic case we shall assume that $\pi$ is an unramified
representation and that $B$ corresponds to the spherical vector.

The generic $\GL_2(F)$ representation $\tau$, however, will be completely arbitrary. In the
non-archimedean case, we will restrict our attention to ramified representations $\tau$,
since the unramified case has been done in \cite{Fu}. In both the archimedean and
non-archimedean cases, $W^\#$ will be a vector in the induced representation which has a
suitable right transformation property under the maximal compact subgroup of $H(F)$
depending on that of the Bessel vector $B$.
\subsection{Review of Furusawa's theory}\label{furusawa-review}
In this section we recall some of the theory of \cite{Fu} relevant for this paper.
For simplicity we work over the ground field $\Q$.
\subsubsection*{Bessel models}
Let $a,b,c\in\Q$ such that $D=4ac-b^2>0$ is a non-square in $\Q^\times$.
Then $L=\Q(\sqrt{-D})$ is an imaginary quadratic field extension of $\Q$. Let $\A_L$ be the ring of adeles of $L$. The adelic points $T(\A)$ of the group defined in (\ref{TFdefeq}) satisfies  $T(\A)\cong\A_L^\times$
and $T(\Q)\simeq L^{\times}$.

We fix the additive character $\psi=\prod_p\psi_p$ of $\Q\backslash\A$ for which
$\psi_p$ has conductor $\Z_p$ for all primes $p$ and for which $\psi_\infty(x)=e^{-2\pi ix}$
for $x\in\R$. Let $U$ be the unipotent radical of the Siegel parabolic subgroup
of $H$. Let $\theta$ be the character of $U(\A)$ given by
$$
 \theta(\mat{1}{X}{}{1})=\psi({\rm tr}(SX)),\qquad
 X\in M_2(\A),\;X=\,^tX.
$$
As in Sect.\ \ref{besselsubgroupsec}, $R=TU$ is the
Bessel subgroup defined by $S$. Let $\Lambda$ be a Hecke character of $L$, i.e.,
a character of $T(\Q)\backslash T(\A)\cong L^\times\backslash\A_L^\times$.
Then the map
$$
 tu\longmapsto \Lambda(t)\theta(u),\qquad t\in T(\A),\:u\in U(\A),
$$
is a character of $R(\Q)\backslash R(\A)$, which we denote by $\Lambda\otimes\theta$.

Let $\pi=\otimes\pi_p$ be a cuspidal, automorphic representation of $H(\A)$.
Let $V_\pi$ be the space of automorphic forms realizing $\pi$.
Assume that a Hecke character $\Lambda$ as above is chosen such that
the restriction of $\Lambda$ to $\A^\times$ coincides with $\omega_\pi$, the
central character of $\pi$.
For each $\phi\in V_\pi$ consider the corresponding Bessel function
\begin{equation}\label{Bphieq}
 B_\phi(g)=\int\limits_{Z_H(\A)R(\Q)\backslash R(\A)}
 (\Lambda\otimes\theta)(r)^{-1}\phi(rg)\,dr,
\end{equation}
where $Z_H$ is the center of $H$. If one of these integrals is non-zero, then all
are non-zero, and we obtain a model $\mathcal{B}_{\Lambda,\theta,\psi}(\pi)$
of $\pi$ consisting of functions on $H(\A)$
with the obvious transformation property on the left with respect to $R(\A)$.
In this case, we say that $\pi$ has a \emph{global Bessel model} of type $(S,\Lambda,\psi)$.
It implies that the \emph{local} Bessel model
$\mathcal{B}_{\Lambda_p,\theta_p,\psi_p}(\pi_p)$ exists for every $p$. In fact, there is a
canonical isomorphism
$$
 \bigotimes\limits_p\mathcal{B}_{\Lambda_p,\theta_p,\psi_p}(\pi_p)\cong
 \mathcal{B}_{\Lambda,\theta,\psi}(\pi).
$$
If $(B_p)_p$ is a collection of local Bessel functions $B_p\in
\mathcal{B}_{\Lambda_p,\theta_p,\psi_p}(\pi_p)$ such that $B_p\big|_{H(\Z_p)}=1$
for almost all $p$, then this isomorphism is such that $\otimes_pB_p$ corresponds
to the global function
\begin{equation}\label{locglobBreleq}
 B(g)=\prod_pB_p(g_p),\qquad g=(g_p)_p\in H(\A).
\end{equation}
\subsubsection*{Global induced representations}
The Eisenstein series $E(h,s)$ entering into the global
integral (\ref{globalintegraleq}) below will be defined from a section in a global
induced representation of $G(\A)$.
We therefore now discuss various models of such induced representations.
Let $(\tau,V_\tau)$ be a cuspidal, automorphic representation of $\GL_2(\A)$.
Let $\chi_0$ be a character of $L^\times\backslash\A_L^\times$ such that
the restriction of $\chi_0$ to $\A^\times$ concides with $\omega_\tau$, the
central character of $\tau$. Then, as in (\ref{M2representationseq}) in the local case,
$\chi_0$ can be used to extend $\tau$ to a representation of $M^{(2)}(\A)$,
denoted by $\chi_0 \times \tau$. Let $\chi$ be another character of
$L^\times\backslash\A_L^\times$, considered as a character of $M^{(1)}(\A)$.
This data defines a family of induced representations
$I(s,\chi,\chi_0,\tau)$ of $G(\A)$ depending on a complex parameter $s$.
The space of $I(s,\chi,\chi_0,\tau)$ consists of functions
$\varphi:\:G(\A)\rightarrow V_\tau$ with the transformation property
$$
 \varphi(m_1m_2ng)=\delta_P(m_1m_2)^{s+1/2}\chi(m_1)(\chi_0 \times \tau)(m_2)\varphi(g),\qquad
 m_1\in M^{(1)}(\A),\:m_2\in M^{(2)}(\A),\:n\in N(\A).
$$
Since the representation $\tau$ is given as a space of automorphic forms, we may realize
$I(s,\chi,\chi_0,\tau)$ as a space of $\C$-valued functions on $G(\A)$.
More precisely, to each $\varphi$ as above we may attach the function
$f_\varphi$ on $G(\A)$ given by $f_\varphi(g)=(\varphi(g))(1)$.
Each function $f_\varphi$ has the property that $\GL_2(\A)\ni h\mapsto f_\varphi(hg)$
is an element of $V_\tau$, for each $g\in G(\A)$. Let $I_\C(s,\chi,\chi_0,\tau)$
be the model of $I(s,\chi,\chi_0,\tau)$ thus obtained.
A third model of the same representation is obtained by attaching to
$f\in I_\C(s,\chi,\chi_0,\tau)$ the function
\begin{equation}\label{Wffreleq}
 W_f(g) = \int\limits_{\Q \backslash \A}
 f\Big(\begin{bmatrix}1&&&\\&1&&x\\&&1&\\&&&1\end{bmatrix}g\Big)\psi(cx)dx,
 \qquad g\in G(\A).
\end{equation}
The map $f\mapsto W_f$ is injective since $\tau$ is cuspidal.
In fact, $f$ can be recovered from $W_f$ via the formula
\begin{equation}\label{fWfreleq}
 f(g)=\sum_{\lambda\in\Q^\times}W_f
 \Big(\begin{bmatrix}1&&&\\&\lambda\\&&\lambda\\&&&1\end{bmatrix}g\Big),
 \qquad g\in G(\A).
\end{equation}
Let $I_W(s,\chi,\chi_0,\tau)$ be the space of all functions $W_f$.
Now write $\tau\cong\otimes\tau_p$ with local representations $\tau_p$ of $\GL_2(\Q_p)$.
We also factor $\chi=\otimes\chi_p$ and $\chi_0=\otimes\chi_{0,p}$, where
$\chi_p$ and $\chi_{0,p}$ are characters of $\prod_{v|p}L_v^\times$.
Then there are isomorphisms
\begin{equation}\label{locglobindrepdiagrameq}
 \begin{CD}
  I(s,\chi,\chi_0,\tau)@>\sim>>\otimes_p I(s,\chi_p,\chi_{0,p},\tau_p)\\
  @V\sim VV @VV=V\\
  I_\C(s,\chi,\chi_0,\tau)@>\sim>>\otimes_p I(s,\chi_p,\chi_{0,p},\tau_p)\\
  @V\sim VV @VV\sim V\\
  I_W(s,\chi,\chi_0,\tau)@>\sim>>\otimes_p I_W(s,\chi_p,\chi_{0,p},\tau_p)
 \end{CD}
\end{equation}
Here, the local induced representation $I(s,\chi_p,\chi_{0,p},\tau_p)$ consists
of functions taking values in a model $V_{\tau_p}$ of $\tau_p$; see Sect.\
\ref{parabolicinductionsec} for the precise definition. Assume that
$V_{\tau_p}=\mathcal{W}(\tau_p,\psi_p^{-c})$ is the Whittaker model of $V_{\tau_p}$
with respect to the additive character $\psi_p^{-c}$. If we attach to each
$f_p\in I(s,\chi_p,\chi_{0,p},\tau_p)$ the function
$W_{f_p}(g)=f_p(g)(1)$, then we obtain the model $I_W(s,\chi_p,\chi_{0,p},\tau_p)$
of the same induced representation. The bottom isomorphism in diagram
(\ref{locglobindrepdiagrameq}) is such that if
$W_p\in I_W(s,\chi_p,\chi_{0,p},\tau_p)$ are given, with the property that
$W_p\big|_{G(\Z_p)}=1$ for almost all $p$, then the corresponding element of
$I_W(s,\chi,\chi_0,\tau)$ is the function
\begin{equation}\label{Wfactorizationeq}
 W(g)=\prod_{p\leq\infty}W_p(g_p),\qquad g=(g_p)_p\in G(\A).
\end{equation}
\subsubsection*{The global integral and the basic identity}
As above let $(\tau,V_\tau)$ be a cuspidal, automorphic representation of $\GL_2(\A)$,
extended to a representation of $M^{(2)}(\A)$ via a character $\chi_0$ of
$L^\times\backslash\A_L^\times$. Let further $(\pi,V_\pi)$
be a cuspidal, automorphic representation of $H(\A)$ which has a global Bessel
model of type $(S,\Lambda,\psi)$, where $S=\mat{a}{b/2}{b/2}{c}$ as above.
Define the character $\chi$ of $L^\times\backslash\A_L^\times$ by
$$
 \chi(a)=\Lambda(\bar a)^{-1}\chi_0(\bar a)^{-1},\qquad a\in\A_L^\times.
$$
Let $f(g,s)$ be an analytic family in $I_\C(s,\chi,\chi_0,\tau)$.
For ${\rm Re}(s)$ large enough we can form the Eisenstein series
$$
 E(g,s;f)=\sum_{\gamma\in P(\Q)\backslash G(\Q)}f(\gamma g,s).
$$
In fact, $E(g,s;f)$ has a meromorphic continuation to the entire plane.
In \cite{Fu} Furusawa studied integrals of the form
\begin{equation}\label{globalintegraleq}
 Z(s,f,\phi)=\int\limits_{H(\Q)Z_H(\A)\backslash H(\A)}E(h,s;f)\phi(h)\,dh,
\end{equation}
where $\phi\in V_\pi$. Theorem (2.4) of \cite{Fu}, the ``Basic Identity'', states that
\begin{equation}\label{basicidentityeq}
 Z(s,f,\phi)=\int\limits_{R(\A)\backslash H(\A)}W_f(\eta h,s)B_\phi(h)\,dh,
 \qquad\eta=\begin{bmatrix}1\\\alpha&1\\&&1&-\bar\alpha\\&&&1\end{bmatrix},
 \qquad\alpha=\frac{b+\sqrt{d}}{2c},
\end{equation}
where $R(\A)$ is the Bessel subgroup determined by $(S,\Lambda,\psi)$, and
$B_\phi$ is the Bessel function corresponding to $\phi$; see (\ref{Bphieq}).
The function $W_f(\,\cdot\,,s)$
appearing in (\ref{basicidentityeq}) is the element of $I_W(s,\chi,\chi_0,\tau)$
corresponding to $f(\,\cdot\,,s)\in I_\C(s,\chi,\chi_0,\tau)$; see (\ref{Wffreleq})
for the formula relating $f$ and $W_f$.
\subsubsection*{Factorization}
The importance of the basic identity lies in the fact that the integral
on the right side of (\ref{basicidentityeq}) is Eulerian.
Namely, assume that $f(\,\cdot\,,s)$
corresponds to a pure tensor $\otimes f_p$ via the middle isomorphism in
(\ref{locglobindrepdiagrameq}). Assume that $W_p\in I_W(s,\chi_p,\chi_{0,p},\tau_p)$
corresponds to $f_p\in I(s,\chi_p,\chi_{0,p},\tau_p)$. Then
$$
 W_f(g,s)=\prod_{\p\leq\infty}W_p(g_p,s),\qquad g=(g_p)_p\in G(\A),
$$
see (\ref{Wfactorizationeq}). Assume further that the global Bessel function $B_\phi$
factorizes as in (\ref{locglobBreleq}). Then it follows from (\ref{basicidentityeq}) that
\begin{equation}\label{Zfactorizationeq}
 Z(s,f,\phi)=\prod_{p\leq\infty}Z_p(s,W_p,B_p),
\end{equation}
where
\begin{equation}\label{localzetaeq}
 Z_p(s,W_p,B_p)=\int\limits_{R(\Q_p)\backslash H(\Q_p)}W_p(\eta h,s)B_p(h)\,dh.
\end{equation}
Furusawa has calculated the local integrals (\ref{localzetaeq}) in the case where all
the data is unramified. The result is that
\begin{equation}\label{unramifiedlocalzetaeq}
 Z_p(s,W_p,B_p)=\frac{L(3s+\frac12,\tilde\pi_p\times\tilde\tau_p)}
  {L(6s+1,\chi_p\big|_{\Q_p^\times})L(3s+1,\tau_p\times\mathcal{AI}(\Lambda_p)
  \times(\chi_p\big|_{\Q_p^\times}))}.
\end{equation}
Here, $\mathcal{AI}(\Lambda_p)$ denotes the representation of $\GL_2(\Q_p)$ obtained
from the character $\Lambda_p$ via automorphic induction. A $\sim$ over a representation
denotes its contragredient.
Therefore, up to finitely many factors and up to $L$-functions with well-known
analytic properties, the global integral (\ref{globalintegraleq}) represents
the $\GSp_4\times\GL_2$ $L$-function $L(s,\tilde\pi\times\tilde\tau)$.
\nl
In order to obtain more detailed information about the analytic properties of
$L(s,\tilde\pi\times\tilde\tau)$, one has to take the ramified and archimedean
places into account as well. As is evident from (\ref{localzetaeq}), the
choice of local vectors $W_p\in I_W(s,\chi_p,\chi_{0,p},\tau_p)$ and
$B_p\in\mathcal{B}_{\Lambda_p,\theta_p,\psi_p}(\pi_p)$ is crucial.
Any time these vectors can be chosen such that all local integrals (\ref{localzetaeq})
are non-zero, one obtains an integral representation for $L(s,\tilde\pi\times\tilde\tau)$,
with possibly finitely many undesirable Euler factors which need to be controlled
as well.
\section{Local non-archimedean theory}\label{non-arch-section}
In this section we evaluate the local zeta integral (\ref{localZseq}) in the
non-archimedean setting. The key steps are the choices of the
vector $W^\#$ and the actual computation of the integral $Z(s,W^\#,B)$. The vector $B$ will be chosen to be the spherical vector in $\pi$. For $W^\#$, we want to choose a vector in the induced representation that is right invariant under $K^H = H(\OF)$. We will show that such vectors exist and will obtain a canonical one using the newform theory for $\GL_2$.
\subsection{Notation}
Let $F$ be a non-archimedean local field of characteristic zero.
Let $\OF$, $\p$, $\varpi$, $q$ be the ring of integers, prime ideal, uniformizer
and cardinality of the residue class field $\OF/\p$, respectively. Recall that we fix three
elements $a,b,c\in F$ such that $d := b^2-4ac \neq 0$. Let $L$ be as in (\ref{Ldefeq}).
We shall make the following {\bf assumptions}:
\begin{description}
 \item[(A1)] $a,b\in\OF$ and $c\in\OF^\times$.
 \item[(A2)] If $d \not\in F^{\times2}$, then $d$ is the generator of the discriminant of $L/F$.
  If $d \in F^{\times2}$, then $d \in \OF^{\times}$.
\end{description}
We set the Legendre symbol as follows,
\begin{equation}\label{legendresymboldefeq}
 \Big(\frac L{\p}\Big) := \left\{
                  \begin{array}{l@{\qquad\text{if }}l@{\qquad}l}
                    -1, &d \not\in F^{\times2},\:d \not\in \p&\mbox{(the inert case)}, \\
                    0, &d\not\in F^{\times2},\:d\in\p&\mbox{(the ramified case)}, \\
                    1, &d \in F^{\times2}&\mbox{(the split case)}.
                  \end{array}\right.
\end{equation}
If $L$ is a field, then let $\OF_L$ be its ring of integers.
If $L = F \oplus F$, then let $\OF_L = \OF \oplus \OF$.
Note that $x\in \OF_L$ if and only if $N(x), \tr(x) \in \OF$.
If $L$ is a field then we have $x \in \OF_L^{\times}$ if and only if
$N(x) \in \OF^{\times}$. If $L$ is not a field then
$x \in \OF_L, N(x) \in \OF^{\times}$ implies that
$x \in \OF_L^{\times} = \OF^{\times} \oplus \OF^{\times}$.
Let $\varpi_L$ be the uniformizer of $\OF_L$ if $L$ is a field and set
$\varpi_L = (\varpi,1)$ if $L$ is not a field.
Note that, if $(\frac{L}{\p}) \neq -1$, then $N(\varpi_L) \in \varpi \OF^{\times}$. Let
$\alpha\in\OF_L$ be as in (\ref{alphadefeq}). Then, by Lemma 3.1.1 of \cite{PS1},
\begin{equation}\label{OL-in-terms-of-alpha}
 \OF_L = \OF + \alpha \OF.
\end{equation}
We fix the following ideal in $\OF_L$,
\begin{equation}\label{ideal defn}\renewcommand{\arraystretch}{1.3}
 \P := \p\OF_L = \left\{
                  \begin{array}{l@{\qquad\text{if }}l}
                    \p_L & \big(\frac L{\p}\big) = -1,\\
                    \p_L^2 & \big(\frac L{\p}\big) = 0,\\
                    \p \oplus \p & \big(\frac L{\p}\big) = 1.
                  \end{array}
                \right.
\end{equation}
Here, $\p_L$ is the maximal ideal of $\OF_L$ when $L$ is a field
extension. Note that $\P$ is prime only if $\big(\frac
L\p\big)=-1$. We have $\P^n\cap\OF=\p^n$ for all $n\geq0$.
\subsection{The spherical Bessel function}\label{sphericalbesselfunctionsec}
Let $(\pi,V_\pi)$ be an unramified, irreducible, admissible
representation of $H(F)$. Then $\pi$ can be realized as the
unramified constituent of an induced representation of the form
$\chi_1\times\chi_2\rtimes\sigma$, where $\chi_1$, $\chi_2$ and
$\sigma$ are unramified characters of $F^\times$; here, we used
the notation of \cite{ST} for parabolic induction. Let
$$
 \gamma^{(1)}=\chi_1\chi_2\sigma,\qquad\gamma^{(2)}=\chi_1\sigma,\qquad
 \gamma^{(3)}=\sigma,\qquad\gamma^{(4)}=\chi_2\sigma.
$$
Then $\gamma^{(1)}\gamma^{(3)}=\gamma^{(2)}\gamma^{(4)}$ is the central character
of $\pi$. The numbers $\gamma^{(1)}(\varpi),\ldots,\gamma^{(4)}(\varpi)$ are the
Satake parameters of $\pi$. The degree-$4$ $L$-factor of $\pi$ is given by
$\prod_{i=1}^4(1-\gamma^{(i)}(\varpi)q^{-s})^{-1}$.

Let $\Lambda$ be any character of $T(F)\cong L^\times$. We assume
that $V_\pi$ is the Bessel model with respect to the character
$\Lambda\otimes\theta$ of $R(F)$; see Sect.\
\ref{besselsubgroupsec}. Let $B\in V_\pi$ be a spherical vector.
By \cite{Su}, Proposition 2-5, we have $B(1)\neq0$, which implies that necessarily
$\Lambda\big|_{\OF_L^\times}=1$. For $l,m\in\Z$ let
\begin{equation}\label{hlmdefeq}
 h(l,m)=\begin{bmatrix}\varpi^{2m+l}\\&\varpi^{m+l}\\&&1\\&&&\varpi^m\end{bmatrix}.
\end{equation}
Then, as in (3.4.2) of \cite{Fu},
\begin{equation}\label{RFKHrepresentativeseq}
 H(F)=\bigsqcup_{l\in\Z}\bigsqcup_{m\geq0}R(F)h(l,m)K^H,\qquad
 K^H=H(\OF).
\end{equation}
The double cosets on the right hand side are pairwise disjoint.
Since $B$ transforms on the left under $R(F)$ by the character
$\Lambda \otimes \theta$ and is right $K^H$-invariant, it follows
that $B$ is determined by the values $B(h(l,m))$. By Lemma (3.4.4)
of \cite{Fu} we have $B(h(l,m))=0$ for $l<0$, so that $B$ is
determined by the values $B(h(l,m))$ for $l,m\geq0$.

In \cite{Su}, 2-4, Sugano has given a formula for $B(h(l,m))$ in
terms of a generating function. It turns out that for our purposes
we only require the values $B(h(l,0))$. In this special case
Sugano's formula reads
\begin{equation}\label{suganox0eq}
 \sum_{l\geq0}B(h(l,0))y^l= \frac{H(y)}{Q(y)},
\end{equation}
where
\begin{align}\label{suganoQeq}
 Q(y)&=\prod_{i=1}^4\big(1-\gamma^{(i)}(\varpi)q^{-3/2}y\big)
\end{align}
and
\begin{equation}\label{suganoHeq}
 H(y) = \left\{\begin{array}{ll}
    1-q^{-4}\Lambda(\varpi)y^2& \hbox{ if }
     \big(\frac L\p\big)=-1,\\[1ex]
    1-q^{-2}\Lambda(\varpi_L)y& \hbox{ if }
     \big(\frac L\p\big)= 0,\\[1ex]
    1-q^{-2}\big(\Lambda(\varpi_L)+\Lambda(\varpi\varpi_L^{-1})\big)y
     +q^{-4}\Lambda(\varpi)y^2& \hbox{ if } \big(\frac L\p\big)= 1.
 \end{array}\right.
\end{equation}
\subsection{Double coset decompositions}\label{doublecosetsec}
Let $K^G=G(F)\cap\GL_4(\OF_L)$, a maximal compact
subgroup. We define the principal congruence subgroups
\begin{equation}\label{princ-cong-defn}
 \Gamma(\P^r):=G(F)\cap\begin{bmatrix}1+\P^r&\P^r&\P^r&\P^r\\
 \P^r&1+\P^r&\P^r&\P^r\\\P^r&\P^r&1+\P^r&\P^r\\\P^r&\P^r&\P^r&1+\P^r\end{bmatrix}
\end{equation}
with $\P$ as in (\ref{ideal defn}).
For $r=0$ we understand that $\Gamma(\P^r)=K^G$.

The main result of this section is the double coset decomposition in
Proposition \ref{disj-doub-coset-decomp-general-n-prop} below. In Proposition \ref{disj-doub-coset-decomp-general-n-prop}, we obtain representatives for $P(F) \backslash G(F) / K^H$ and $P(F) \backslash G(F) / K^H \Gamma(\P^r)$, and the corresponding double cosets in $K^G$. This will be crucial for the definition of $W^\#$. For this, the first step is $r=1$. We have to treat the case $\big(\frac{L}{\p}\big) = 1$ separately, for which we need Lemma \ref{gl4finitelemma}. We then obtain the disjointness of various double cosets in Lemma \ref{etamdisjointnesslemma}. Using these, we obtain the $r=1$ case in Lemma \ref{double-coset-decomp-n=1-lemma}. And, finally, the general case is done in Proposition \ref{disj-doub-coset-decomp-general-n-prop}.

We start with
the following lemma, which will be used for the split case of
Lemma \ref{double-coset-decomp-n=1-lemma}.
\begin{lemma}\label{gl4finitelemma}
 Let $k$ be a field. Let $P_4(k)$ be the parabolic subgroup of $\GL_4(k)$ consisting
 of matrices of the form
 $$
  \begin{bmatrix} *&*&*&*\\&*&*&*\\&&*\\&*&*&*\end{bmatrix}.
 $$
 Let $\GSp_4(k)$ be defined using the symplectic form $\mat{0}{1_2}{-1_2}{0}$, considered
 as a subgroup of $\GL_4(k)$. Then
 $$
  \GL_4(k)=P_4(k)\GSp_4(k)\sqcup P_4(k)t_1\GSp_4(k),\qquad\text{where }
  t_1=\begin{bmatrix}&1\\1\\&&1\\&&&1\end{bmatrix}.
 $$
\end{lemma}
{\bf Proof.} It is easy to see that the double cosets represented by $1$ and $t_1$
are disjoint. We have to show that every element of $\GL_4(k)$ lies in one of these
two double cosets. Let $W_4$ be the Weyl group of $\GL_4$. Representatives for the
generators of $W_4$ are given by
\begin{equation}\label{gl4finitelemmaeq1}
 t_1=\begin{bmatrix}&1\\1\\&&1\\&&&1\end{bmatrix},\qquad
 t_2=\begin{bmatrix}1\\&&&1\\&&1\\&-1\end{bmatrix},\qquad
 t_3=\begin{bmatrix}1\\&1\\&&&1\\&&1\end{bmatrix}.
\end{equation}
Note that $t_2\in P_4(k)$. The Bruhat decomposition for $\GL_4(k)$ implies that
\begin{align}\label{gl4finitelemmaeq2}
 \GL_4(k)&=P_4\sqcup P_4t_1\begin{bmatrix}1&*\\&1\\&&1\\&&&1\end{bmatrix}
  \sqcup P_4t_3\begin{bmatrix}1\\&1\\&&1\\&&*&1\end{bmatrix}\nonumber\\
 &\quad\sqcup P_4t_1t_2\begin{bmatrix}1&&&*\\&1&&*\\&&1\\&&&1\end{bmatrix}
  \sqcup P_4t_1t_3\begin{bmatrix}1&*\\&1\\&&1\\&&*&1\end{bmatrix}
  \sqcup P_4t_3t_2\begin{bmatrix}1\\&1&*&*\\&&1\\&&&1\end{bmatrix}\nonumber\\
 &\quad\sqcup P_4t_1t_2t_3\begin{bmatrix}1&&*\\&1&*\\&&1\\&&*&1\end{bmatrix}
  \sqcup P_4t_1t_3t_2\begin{bmatrix}1&&&*\\&1&*&*\\&&1\\&&&1\end{bmatrix}
  \sqcup P_4t_3t_2t_1\begin{bmatrix}1&*&*&*\\&1\\&&1\\&&&1\end{bmatrix}\nonumber\\
 &\quad\sqcup P_4t_1t_2t_3t_2\begin{bmatrix}1&&*\\&1&*&*\\&&1\\&&*&1\end{bmatrix}
  \sqcup P_4t_1t_3t_2t_1\begin{bmatrix}1&*&*&*\\&1&&*\\&&1\\&&&1\end{bmatrix}
  \sqcup P_4t_1t_2t_3t_2t_1\begin{bmatrix}1&*&*&*\\&1&*\\&&1\\&&*&1\end{bmatrix},
\end{align}
where we simply wrote $P_4$ for $P_4(k)$. Each of these cosets can be reduced to one of
the first two by multiplying with suitable elements of $\GSp_4(k)$ on the right. For example,
\begin{align*}
 P_4t_3\begin{bmatrix}1\\&1\\&&1\\&&*&1\end{bmatrix}\GSp_4(k)
  &=P_4t_3\begin{bmatrix}1&*\\&1\\&&1\\&&&1\end{bmatrix}\GSp_4(k)
  =P_4t_3\GSp_4(k)\\
  &=P_4t_1t_1t_3\GSp_4(k)
  =P_4t_1\GSp_4(k).
\end{align*}
Similary,
\begin{align}\label{gl4finitelemmaeq3}
 P_4t_3t_2t_1\begin{bmatrix}1&*&*&*\\&1\\&&1\\&&&1\end{bmatrix}\GSp_4(k)
  &=P_4t_2t_3t_2t_1\begin{bmatrix}1&*&&*\\&1\\&&1\\&&&1\end{bmatrix}\GSp_4(k)
  =P_4t_3t_2t_3t_1\begin{bmatrix}1&*&&*\\&1\\&&1\\&&&1\end{bmatrix}\GSp_4(k)\nonumber\\
  &=P_4t_3t_2\begin{bmatrix}1\\ *&1&*\\&&1&\\&&&1\end{bmatrix}\GSp_4(k) =P_4t_3\begin{bmatrix}1\\&1\\&&1&\\ *&&*&1\end{bmatrix}\GSp_4(k).
\end{align}
If the element in the lower left corner is zero, then we are reduced to the case worked out above. If the element in the lower left corner of (\ref{gl4finitelemmaeq3}) is non-zero,
then, using
\begin{equation}\label{gl4finitelemmaeq4}
 \mat{1}{}{x}{1}=\mat{-x^{-1}}{1}{}{x}\mat{}{1}{1}{}\mat{1}{x^{-1}}{}{1},
\end{equation}
we get
\begin{align*}
 P_4t_3\begin{bmatrix}1\\&1\\&&1&\\ *&&*&1\end{bmatrix}\GSp_4(k)
  &=P_4t_3\begin{bmatrix} *&&&*\\&1\\&&1&\\&&&*\end{bmatrix}
   t_2t_1t_2\begin{bmatrix}1&&&*\\&1\\&&1&\\&&&1\end{bmatrix}
   \begin{bmatrix}1\\&1\\&&1&\\&&*&1\end{bmatrix}\GSp_4(k)\\
  &=P_4t_3t_2t_1t_2\begin{bmatrix}1&&&*\\&1\\&&1&\\&&&1\end{bmatrix}
   \begin{bmatrix}1\\&1\\&&1&\\&&*&1\end{bmatrix}\GSp_4(k)\\
  &=P_4t_3t_2t_1t_2\begin{bmatrix}1&&&*\\&1\\&&1&\\&&&1\end{bmatrix}\GSp_4(k)
  =P_4t_3t_2t_1t_2\begin{bmatrix}1\\&1&*\\&&1&\\&&&1\end{bmatrix}\GSp_4(k)\\
  &=P_4t_3t_2t_1t_2\GSp_4(k)
  =P_4t_2t_3t_2t_1\GSp_4(k)
  =P_4t_3t_2t_3t_1\GSp_4(k)\\
  &=P_4t_3\GSp_4(k)
  =P_4t_1\GSp_4(k).
\end{align*}
In a similar way, all cosets occuring in (\ref{gl4finitelemmaeq2}) can be reduced to
one of the first two after multiplication on the right with $\GSp_4(k)$.
Finally, for elements of the second coset, we have
$$
 P_4t_1\begin{bmatrix}1&*\\&1\\&&1\\&&&1\end{bmatrix}\GSp_4(k)
  =P_4t_1\begin{bmatrix}1\\&1\\&&1\\&&*&1\end{bmatrix}\GSp_4(k)
  =P_4t_1\GSp_4(k).
$$
This concludes the proof.\qed

We return to the group $G(F)$. Recall that
\begin{equation}\label{etadef2eq}
 \eta = \begin{bmatrix}1&&&\\ \alpha&1&&\\&&1&-\bar{\alpha}\\&&&1\end{bmatrix},
 \qquad\alpha\text{ as in }(\ref{alphadefeq}).
\end{equation}
For any $m \geq 0$, we let
\begin{equation}\label{etamdef2eq}
 \eta_m = \begin{bmatrix}1&0&&\\\alpha \varpi^m&1&&\\
  &&1&-\bar{\alpha}\varpi^m\\&&0&1\end{bmatrix}.
\end{equation}
For systematic reasons, we let $\eta_\infty$ be the identity matrix.
\begin{lemma}\label{etamdisjointnesslemma}
 Let $K^G=G(F)\cap\GL_4(\OF_L)$ as before.
 \begin{enumerate}
  \item The subsets of $K^G$ given by
   \begin{equation}\label{etamdisjointnesslemmaeq0}
    P(\OF)\eta_mK^H,\qquad m\in\{0,1,2,\ldots,\infty\},
   \end{equation}
   are pairwise disjoint.
  \item Let $r\geq1$. The subsets of $K^G$ given by
   \begin{equation}\label{etamdisjointnesslemmaeq0b}
    P(\OF)\eta_mK^H\Gamma(\P^r),\qquad m\in\{0,\ldots,r\},
   \end{equation}
   are pairwise disjoint.
 \end{enumerate}
\end{lemma}
{\bf Proof.} i) Let $g=p\eta_m k$ with $k\in K^H$ and
$$
 p=\begin{bmatrix}\zeta\\&a&&b\\&&\bar\zeta^{-1}\mu\\&c&&d\end{bmatrix}
 \begin{bmatrix}1&z\\&1\\&&1\\&&-\bar z&1\end{bmatrix}
 \begin{bmatrix}1&&x&y\\&1&\bar y\\&&1\\&&&1\end{bmatrix}\in P(\OF)
 \qquad(\mu=\bar ad-b\bar c).
$$
A calculation shows that the $(3,2)$-coefficient of $gJ\,^t\!g$ is given by
\begin{equation}\label{etamdisjointnesslemmaeq1}
 (gJ\,^t\!g)_{3,2}=a\mu\varpi^m\bar\zeta^{-1}(\bar\alpha-\alpha)
\end{equation}
and that the $(3,4)$-coefficient of $gJ\,^t\!g$ is given by
\begin{equation}\label{etamdisjointnesslemmaeq2}
 (gJ\,^t\!g)_{3,4}=c\mu\varpi^m\bar\zeta^{-1}(\bar\alpha-\alpha)
\end{equation}
(with the understanding that the right sides of
(\ref{etamdisjointnesslemmaeq1}) and (\ref{etamdisjointnesslemmaeq2}) are zero
if $m=\infty$). Since at least one of $a$ or $c$ is in $\OF_L^\times$, it follows
that the function on $K^G$ given by
\begin{equation}\label{etamdisjointnesslemmaeq4}
 g\longmapsto{\rm min}\big(v((gJ\,^t\!g)_{3,2}),v((gJ\,^t\!g)_{3,4})\big)
\end{equation}
takes different values on the double cosets (\ref{etamdisjointnesslemmaeq0}).
\nll
ii) The argument is similar as in i); one considers the valuation of the
$(3,2)$- and the $(3,4)$-coefficient of $gJ\,^t\!g$ mod $\P^r$.
\qed

\begin{lemma}\label{princcongsubgrpdecomplemma}
 Let $r$ be a positive integer and
 \begin{equation}\label{princcongsubgrpdecomplemmaeq1}
  \gamma\in G(F)\cap\begin{bmatrix}\OF_L^\times&\P^r&\OF_L&\OF_L\\
  \OF_L&\OF_L^\times&\OF_L&\OF_L\\\P^r&\P^r&\OF_L^\times&\OF_L\\\P^r&\P^r&\P^r&\OF_L^\times
  \end{bmatrix}.
 \end{equation}
 Then $\gamma$ can be written in the form
 $$
  \gamma=p\eta_m h,\qquad p\in P(\OF),\;h\in K^H,\;\eta_m
  =\begin{bmatrix}1\\\alpha\varpi^m&1\\&&1&-\bar\alpha\varpi^m\\&&&1\end{bmatrix},
 $$
 with a uniquely determined $m\in\{0,1,2,\ldots,\infty\}$.
\end{lemma}
{\bf Proof.} The uniqueness of $m$ follows from Lemma \ref{etamdisjointnesslemma}; we
will show that such an $m$ exists. The group (\ref{princcongsubgrpdecomplemmaeq1})
has an Iwahori decomposition, enabling us to write
$$
 \gamma=\mat{1}{B}{}{1}\mat{A}{}{}{^t\!\bar A^{-1}}\mat{1}{}{C}{1},\qquad
 A\in\mat{\OF_L^\times}{\P^r}{\OF_L}{\OF_L^\times},\;
 B\in\mat{\OF_L}{\OF_L}{\OF_L}{\OF_L},\;C\in\mat{\P^r}{\P^r}{\P^r}{\P^r}.
$$
Decomposing $A$ further in the form $\mat{\OF_L^\times}{\P^r}{}{\OF_L^\times}
\mat{1}{}{\OF_L}{1}$, and multiplying on the left with an appropriate element of $P(\OF)$, we may assume that
$$
 \gamma=\begin{bmatrix}1\\z&1\\&&1&-\bar z\\&&&1\end{bmatrix}
 \begin{bmatrix}1\\&1\\x_1&\bar y&1\\y&x_2&&1\end{bmatrix},\qquad
 x_1,x_2\in\P^r\cap F,\;y\in\P^r,\;z\in\OF_L.
$$
Observing that $\OF_L=\OF+\alpha\OF$ and multiplying on the right with an appropriate
elements of $K^H$, we may assume that
\begin{equation}\label{princcongsubgrpdecomplemmaeq2}
 \gamma=\begin{bmatrix}1\\\alpha z&1\\&&1&-\bar\alpha z\\&&&1\end{bmatrix}
 \begin{bmatrix}1\\&1\\&\bar\alpha y&1\\\alpha y&&&1\end{bmatrix},\qquad y,z\in\OF.
\end{equation}
If $y=z=0$, then $\gamma=1=\eta_\infty$, and we are done. Assume that $y$ and $z$
are not both zero. Then the identities
\begin{equation}\label{princcongsubgrpdecomplemmaeq3}
 \begin{bmatrix}1\\&1\\&&1\\&-yz^{-1}&&1\end{bmatrix}
 \begin{bmatrix}1\\\alpha z&1\\&&1&-\bar\alpha z\\&&&1\end{bmatrix}
 \begin{bmatrix}1\\&1\\&\bar\alpha y&1\\\alpha y&&&1\end{bmatrix}
 =\begin{bmatrix}1\\\alpha z&1\\&&1&-\bar\alpha z\\&&&1\end{bmatrix}
 \begin{bmatrix}1\\&1\\-\alpha\bar\alpha yz&&1\\&-yz^{-1}&&1\end{bmatrix}
\end{equation}
(for $yz^{-1}\in\OF$)
and
\begin{equation}\label{princcongsubgrpdecomplemmaeq4}
 \begin{bmatrix}1\\&&&1\\&&-1\\&1&&-zy^{-1}\end{bmatrix}\!
 \begin{bmatrix}1\\\alpha z&1\\&&1&-\bar\alpha z\\&&&1\end{bmatrix}\!
 \begin{bmatrix}1\\&1\\&\bar\alpha y&1\\\alpha y&&&1\end{bmatrix}
 =\begin{bmatrix}1\\\alpha y&1\\&&1&-\bar\alpha y\\&&&1\end{bmatrix}
 \begin{bmatrix}1\\&&&1\\\alpha\bar\alpha yz&&-1\\&1&&-zy^{-1}\end{bmatrix}
\end{equation}
(for $zy^{-1}\in\OF$) show that we may assume
\begin{equation}\label{princcongsubgrpdecomplemmaeq5}
 \gamma=\begin{bmatrix}1\\\alpha z&1\\&&1&-\bar\alpha z\\&&&1\end{bmatrix},
 \qquad z\in\OF.
\end{equation}
By using appropriate unit diagonal matrices, we see that such a $\gamma$ defines
the same element of the double coset space
$P(F)\backslash G/K^H$ as $\eta_m$, where $m=v(z)$.\qed

\begin{lemma}\label{double-coset-decomp-n=1-lemma}
 We have the disjoint union
 \begin{equation}\label{double-coset-decomp-n=1-eqn}
  K^G = P(\OF)K^H\Gamma(\P)\,\sqcup\,P(\OF)\eta K^H\Gamma(\P).
 \end{equation}
\end{lemma}
{\bf Proof.} The disjointness follows from Lemma \ref{etamdisjointnesslemma} ii).
To prove that each element of $K^G$ is contained in one of the double cosets,
we will distinguish three cases depending on the value of $\big(\frac L{\p}\big)$.

Let us first assume that $\big(\frac L{\p}\big) = -1$. In this case, $\P$ is the maximal
ideal in $\OF_L$.  Using $K^G/\Gamma(\P) \simeq G(\OF/\p)$ and the Bruhat
decomposition of the group $G(\OF/\p)$, we see that a set $S$ of representatives
of $P(\OF) \backslash K^G /K^H\Gamma(\P)$ can be chosen from
$\{w b : w \in W,\:b \in B(\OF)\}$.
Here, $W$ is the eight element Weyl group of $\GU(2,2)$ and $B$ is the Borel subgroup.
Since $wb = (w b w^{-1}) w$ and $w \in K^H$, we see that $S$ can be chosen from the
opposite of the Borel subgroup. Since diagonal elements are in $P(\OF)$,
elements of $S$ can be chosen of the form
$$
 \begin{bmatrix}1&&&\\z&1&&\\&&1&-\bar{z}\\&&&1\end{bmatrix}
 \begin{bmatrix}1&&&\\&1&&\\x_1&\bar{y}&1&\\y&x_2&&1\end{bmatrix},\qquad
  z, y \in \OF_L,\;x_1, x_2 \in \OF.
$$
Since $\OF_L = \OF + \alpha \OF$ and we can modify the elements of $S$
by elements of $K^H$ on the right, it follows that the elements of $S$ can be
chosen of the form
$$
 \begin{bmatrix}1&&&\\ \alpha z&1&&\\&&1&-\bar{\alpha}z\\&&&1\end{bmatrix}
 \begin{bmatrix}1&&&\\&1&&\\&\bar{\alpha} y&1&\\ \alpha y&&&1\end{bmatrix},\qquad
 z, y \in \OF.
$$
Using (\ref{princcongsubgrpdecomplemmaeq3}) and (\ref{princcongsubgrpdecomplemmaeq4}),
we see that the elements of $S$ can be chosen of the form
$$
 \begin{bmatrix}1&&&\\ \alpha z&1&&\\&&1&-\bar{\alpha}z\\&&&1\end{bmatrix},\qquad
 z\in \OF.
$$
Finally, using unit diagonal matrices, we may assume that $z\in\{0,1\}$.
The assertion follows.

Next, let us assume that $\big(\frac L{\p}\big) = 0$. In this case $\P = \p_L^2$, where $\p_L$ is the maximal ideal of $\OF_L$. We also have $\OF_L/\p_L\cong\OF/\p$, and thus
$\OF_L = \OF + \p_L$. Moreover, $K^G/\Gamma(\p_L)\cong K^H/\Gamma(\p)$, so that
$$
 K^G=K^H \Gamma(\p_L).
$$
The coset representatives for $\Gamma(\p_L) / \Gamma(\P)$ are given by matrices in $K^G$,
where the diagonal entries are in $1+\p_L$ and the off-diagonal entries are in $\p_L$.
It is easy to show that any matrix $g$ in $\Gamma(\p_L)$ can be written as a product
$g_1 g_2 g_3$, where $g_1 \in U(F) \cap \Gamma(\p_L)$,
$g_3 \in \bar{U}(F) \cap \Gamma(\p_L)$ and $g_2$ is a diagonal matrix in $\Gamma(\p_L)$.
Here, $U(F)$ is the unipotent radical of the Borel subgroup and $\bar{U}(F)$ is the
opposite of $U(F)$. Since $U(F) \cap \Gamma(\p_L)$ and diagonal matrices are
contained in $P(\OF)$, and
since we can modify coset representatives of $P(\OF)\backslash K^G/K^H \Gamma(\P)$
by elements of $K^H$ on the right, such a set $S$ of representatives can be chosen from
$$
 \begin{bmatrix}1&&&\\ z&1&&\\&&1&-\bar{z}\\&&&1\end{bmatrix}
 \begin{bmatrix}1&&&\\&1&&\\&\bar{y}&1&\\y&&&1\end{bmatrix},\qquad z, y
 \text{ from a set of representatives for }\p_L/\P.
$$
Let $w_0\in\OF$ be the mod $\p$ unique element such that $\alpha - w_0 \in \p_L$.
Note that $\alpha - w_0 \not\in \P$ by Lemma 3.1.1 (ii) of \cite{PS1}.
Hence, we can take the set $\{(\alpha - w_0) x : x \in \OF/\p\}$ as representatives of $\p_L/\P$. Since we can modify elements of $S$ by elements of $K^H$ on the right, we see that $S$ can be chosen from
$$
 \begin{bmatrix}1&&&\\ \alpha z&1&&\\&&1&-\bar{\alpha}z\\&&&1\end{bmatrix}
 \begin{bmatrix}1&&&\\&1&&\\&\bar{\alpha} y&1&\\ \alpha y&&&1\end{bmatrix},\qquad
 z,y\text{ from a set of representatives for }\OF/\p.
$$
The assertion now follows by imitating the steps in the proof of the case
$\big(\frac L{\p}\big) = -1$.

Finally, assume that $\big(\frac L{\p}\big) = 1$. In this case $L=F\oplus F$.
Accordingly, we can write every element $g\in G$ as a pair $(g_1,g_2)$ with
matrices $g_1,g_2\in\GL_4(F)$. The condition $^t\bar gJg=\mu(g)J$ translates into
$g_2=\mu(g)J^{-1}\,^t\!g_1^{-1}J$. Hence, we obtain an isomorphism
\begin{align}\label{double-coset-decomp-n=1-lemmaeq3}
 G&\stackrel{\sim}{\longrightarrow}\GL_4(F)\times\GL_1(F),\nonumber\\
 g=(g_1,g_2)&\longmapsto(g_1,\mu(g)).
\end{align}
Under this isomorphism, the parabolic subgroup $P(F)$ is mapped onto $P_4(F)\times\GL_1(F)$,
where
$$
 P_4(F)=\begin{bmatrix} *&*&*&*\\&*&*&*\\&&*\\&*&*&*\end{bmatrix}\subset\GL_4(F)
$$
is the parabolic subgroup of $\GL_4(F)$ of the same shape as $P$. The group $K^G$
is mapped onto $\GL_4(\OF)\times\OF^\times$. The principal
congruence subgroup $\Gamma(\P)$ is mapped onto $\Gamma_4(\p)\times(1+\p)$, where
$\Gamma_4(\p)$ is the principal congruence subgroup
of level $\p$ in $\GL_4(F)$. And the group $K^H=\GSp_4(\OF)$ is mapped onto
$$
 K^H_4:=\{(g,\mu(g))\in\GL_4(F)\times\GL_1(F):\:g\in\GSp_4(\OF)\}.
$$
Therefore,
\begin{align*}
 P(\OF)\backslash K^G/K^H\Gamma(\P)&\cong (P_4(\OF)\times\OF^\times)
  \backslash(\GL_4(\OF)\times\OF^\times)/K^H_4(\Gamma_4(\p)\times(1+\p))\\
 &\cong P_4(\OF)\backslash\GL_4(\OF)/\GSp_4(\OF)\Gamma_4(\p)\\
 &\cong P_4(\OF/\p)\backslash\GL_4(\OF/\p)/\GSp_4(\OF/\p).
\end{align*}
By Lemma \ref{gl4finitelemma}, this double coset space is represented by the
elements $1$ and $t_1$. The assertion therefore follows from the easily checked
fact that the element $\eta\in K^G$ maps to an element representing the same
double coset as $t_1$ in $P_4(\OF/\p)\backslash\GL_4(\OF/\p)/\GSp_4(\OF/\p)$.
This completes the proof.\qed

\begin{proposition}\label{disj-doub-coset-decomp-general-n-prop}
 Let $\eta_m$ be as in (\ref{etamdef2eq}). Let $\eta_\infty$ be the
 identity matrix. We have the following disjoint double coset decompositions.
 \begin{enumerate}
  \item
   $$
    K^G=\bigsqcup_{0\leq m\leq\infty} P(\OF)\eta_mK^H.
   $$
  \item For any $r\geq0$,
   $$
    K^G=\bigsqcup\limits_{0 \leq m \leq r}P(\OF)\eta_m K^H \Gamma(\P^r).
   $$
  \item
   $$
    G(F) = \bigsqcup\limits_{0 \leq m \leq \infty}P(F)\eta_m K^H.
   $$
  \item For any $r\geq0$,
   $$
    G(F) = \bigsqcup\limits_{0 \leq m \leq r}P(F)\eta_m K^H \Gamma(\P^r).
   $$
 \end{enumerate}
\end{proposition}
{\bf Proof.} Using the Iwasawa decomposition, iii) follows from i) and
iv) follows from ii). In view of the disjointness stated in Lemma
\ref{etamdisjointnesslemma}, ii) follows from i) by multiplying on the right
with $\Gamma(\P^r)$. Hence it is enough to prove i). The disjointness was
already proved in Lemma \ref{etamdisjointnesslemma}.

Let $g\in K^G$. Then, by Lemma \ref{double-coset-decomp-n=1-lemma}, either
$g=p\gamma k$ or $g=p\eta\gamma k$, where $p\in P(\OF)$, $\gamma\in\Gamma(\P)$,
and $k\in K^H$. In the first case we write $\gamma$ according to Lemma
\ref{princcongsubgrpdecomplemma}. In the second case we write $\eta\gamma$
according to Lemma \ref{princcongsubgrpdecomplemma}. The assertion follows.\qed

The following lemma shows that the first $r-1$ double cosets
occuring in i) are the same as those occurring in  ii) of Proposition \ref{disj-doub-coset-decomp-general-n-prop}.
\begin{lemma}\label{KHGammaKHlemma}
 For any $0\leq m<r$, we have
 $$
  P(\OF)\eta_m K^H\Gamma(\P^r)=P(\OF)\eta_m K^H
 $$
 and
 $$
  P(F)\eta_m K^H\Gamma(\P^r)=P(F)\eta_m K^H.
 $$
\end{lemma}
{\bf Proof.} We will prove the second equality; the argument for the first one
is the same. By Proposition \ref{disj-doub-coset-decomp-general-n-prop},
for any $r>0$,
$$
 G(F) = \bigsqcup\limits_{0 \leq m<r}P(F)\eta_m K^H\sqcup X,\qquad
 X=\bigsqcup\limits_{r \leq m\leq\infty}P(F)\eta_m K^H,
$$
and also
$$
 G(F)=\bigsqcup_{0\leq m<r}P(F)\eta_mK^H\Gamma(\P^r)\;\sqcup Y,
 \qquad Y=P(F)\;\eta_rK^H\Gamma(\P^r).
$$
For $m\geq r$, we have $\eta_m\in P(F)K^H\Gamma(\P^r)=P(F)\eta_rK^H\Gamma(\P^r)$.
Hence $X\subset Y$.
Evidently, for $m<r$, we have $P(F)\eta_m K^H\subset P(F)\eta_mK^H\Gamma(\P^r)$.
It follows that $P(F)\eta_m K^H= P(F)\eta_mK^H\Gamma(\P^r)$.\qed
\subsection{New- and oldforms for $\GU(2,2)$}\label{Wsharp-defn-nonarch-sec}
Using the double coset decompositions from the previous section,
we shall determine the structure of the spaces of vectors in the
induced representations $I(s, \chi, \chi_0, \tau)$ invariant under the
groups $K^H\Gamma(\P^r)$, $r\geq0$. Here, $\tau$ is a representation of $\GL_2(F)$,
and $\chi$, $\chi_0$ are appropriately chosen characters of $L^\times$; see
Sect.\ \ref{parabolicinductionsec}. It turns out that these spaces
of invariant vectors are zero if $r<n$, where $\p^n$ is the conductor of $\tau$.
If $r=n$, then the space of invariant vectors is one-dimensional; in this sense
there is a unique newform. For $r>n$, the dimensions of the spaces of invariant
vectors grow quadratically. We start by recalling some familiar $\GL_2$ theory.
\subsubsection*{The $\GL_2$ newform}
We define congruence subgroups of $\GL_2(F)$, as follows. For
$n=0$, let $K^{(0)}(\p^0)=\GL_2(\OF)$. For $n>0$, let
\begin{equation}\label{K'defeq}
 K^{(0)}(\p^n)=\mat{1+\p^n}{\OF}{\p^n}{\OF^\times}.
\end{equation}
The following result is well known (see \cite{Cas}, \cite{De}).
\begin{theorem}\label{GL2newformtheorem}
 Let $(\tau,V_\tau)$ be a generic, irreducible, admissible representation of $\GL_2(F)$.
 Then the spaces
 $$
  V_\tau(n)=\{v\in V_\tau:\:\tau(g)v=v\text{ for all }g\in K^{(0)}(\p^n)\}
 $$
 are non-zero for $n$ large enough. If $n$ is minimal with $V_\tau(n)\neq0$, then
 $\dim(V_\tau(n))=1$. For $r \geq n$, we have $\dim(V_\tau(r)) = r-n+1$.
\end{theorem}
If $n$ is minimal such that $V_\tau(n)\neq0$, then $\p^n$ is called the
\emph{conductor} of $\tau$, and any non-zero vector in $V_\tau(n)$ is called
a \emph{local newform}.
\begin{lemma}\label{GL2-newform-Kirillov}
 Let $(\tau,V_\tau)$ be a generic, irreducible, admissible representation
 of $\GL_2(F)$ with conductor $\p^n$.  We assume that $V_\tau$ is the
 Whittaker model of $\tau$ with respect to the character of $F$
 given by $\psi^{-c}(x)=\psi(-cx)$, where $c\in\OF^\times$. Let $W^{(0)}$ be a
 local newform. Then $W^{(0)}(1)\neq0$. If $W^{(0)}$ is normalized such that
 $W^{(0)}(1)= 1$, then the following formulas hold.
 \begin{enumerate}
  \item If $\tau$ is a supercuspidal representation, or $\tau = \Omega\St_{\GL(2)}$ is a
   twist of the Steinberg representation with a ramified character $\Omega$, or $\tau$ is
   a principal series representation $\alpha\times\beta$ with
   two ramified characters $\alpha, \beta$
   (such that $\alpha \beta^{-1} \neq |\,|^{\pm1}$), then
   $$
    W^{(0)}(\mat{\varpi^l}{}{}{1}) = \left\{\begin{array}{ll}
     1& \hbox{ if } l = 0, \\[1ex]
     0 & \hbox{ if } l \neq 0.\end{array}\right.
   $$
  \item If $\tau=\alpha\times\beta$ is a principal series representation with
   an unramified character $\alpha$ and a ramified character $\beta$, then
   $$
    W^{(0)}(\mat{\varpi^l}{}{}{1}) = \left\{\begin{array}{ll}
     (\beta(\varpi)q^{-1/2})^l& \hbox{ if } l \geq 0, \\[1ex]
     0 & \hbox{ if } l < 0.\end{array}\right.
   $$
  \item If $\tau = \Omega\St_{\GL(2)}$ is a
   twist of the Steinberg representation with an unramified character $\Omega$, then
   $$
    W^{(0)}(\mat{\varpi^l}{}{}{1}) = \left\{\begin{array}{ll}
     (\Omega(\varpi)q^{-1})^l& \hbox{ if } l \geq 0, \\[1ex]
     0 & \hbox{ if } l < 0.\end{array}\right.
   $$
  \item If $\tau=\alpha\times\beta$ is a principal series representation with
   unramified characters $\alpha$ and $\beta$, then
   $$
    W^{(0)}(\mat{\varpi^l}{}{}{1}) = \left\{\begin{array}{ll}
     \displaystyle q^{-l/2}\sum_{k=0}^l\alpha(\varpi)^k\beta(\varpi)^{l-k}
      & \hbox{ if } l \geq 0, \\[1ex]
     0 & \hbox{ if } l < 0.\end{array}\right.
   $$
 \end{enumerate}
\end{lemma}
To prove this lemma, one can use formulas for the local newform with respect to the
congruence subgroup $\GL_2(\OF)\cap\mat{\OF}{\OF}{\p^n}{1+\p^n}$ (given, amongst
other places, in \cite{Sc}), together with the local functional equation.
\subsubsection*{An auxiliary lemma}
We will derive a lemma which will be used in the proof of Theorem \ref{unique-W-theorem}
further below.
\begin{lemma}\label{something-in-pn}
 Let $\alpha$ be as in (\ref{alphadefeq}).
 Let $x \in \OF_L$ be such that $x \in \OF + \P^n$ and $\alpha x\in \OF + \P^n$
 for a non-negative integer $n$. Then $x \in \P^n$.
\end{lemma}
{\bf Proof.} Using (\ref{OL-in-terms-of-alpha}) and (\ref{ideal defn}), first note that $\OF + \P^n = \OF + \alpha \p^n$. Let $x = y+\alpha z$, with $y
\in \OF$ and $z \in \p^n$. Since $\alpha^2 = \alpha b/c - a/c$, we get
$\alpha x = -az/c + \alpha (y+bz/c)$. Now, $\alpha x \in \OF +
\P^n$ implies that $y+bz/c \in \p^n$ and hence, $y \in \p^n$. This
proves that $x \in \P^n$. \qed
\nll
Recall that $K^G$ is the maximal compact subgroup of $G(F)$ and that $K^H=\GSp_4(\OF)$.
Let the principal congruence subgroups $\Gamma(\P^r)$
of $K^G$ be defined as in (\ref{princ-cong-defn}).
For $m \geq 0$, let $\eta_m$ be as in (\ref{etamdef2eq}).
\begin{lemma}\label{well-defined-sub-lemma}
 Let
 $$
  \hat{m} = \begin{bmatrix}\zeta\\&a'&&b'\\&&\mu\bar{\zeta}^{-1}\\&c'&&d'\end{bmatrix} \in M(F)
  \qquad\mbox{and}\qquad \hat{n} = \begin{bmatrix}
                 1 & z &  &  \\
                  & 1 &  &  \\
                  &  & 1 &  \\
                  &  & -\overline{z} & 1 \\
               \end{bmatrix}
               \begin{bmatrix}
                 1 &  & w & y \\
                  & 1 & \overline{y} &  \\
                  &  & 1 &  \\
                  &  &  & 1 \\
               \end{bmatrix} \in N(F).
 $$
 Let $m, r$ be integers such that $r > m \geq 0$. If
 $A := \eta_{m}^{-1} \hat{m}\hat{n} \eta_{m} \in K^H \Gamma(\P^{r})$ then $c' \in \P^{r-m}$ and $a' \bar{\zeta}^{-1}  \in 1+\P^{r-m}$.
\end{lemma}
{\bf Proof.} Suppose $A := \eta_{m}^{-1} \hat{m}\hat{n} \eta_{m} \in K^H \Gamma(\P^{r})$. First note that $K^H \Gamma(\P^{r}) \subset
M_4(\OF+\P^{r})$. Looking at the $(3,2), (4,2)$ coefficients
of $A$, we see that $c', \bar\alpha c'\varpi^{m} \in
\OF+\P^{r}$. By Lemma \ref{something-in-pn}, we obtain $c' \varpi^{m}
\in \P^{r}$ and hence $c' \in \P^{r-m}$, as required.

Note that $\hat{m} \hat{n} \in K^G$ and $c' \in \P^{r-m} \subset \P$
implies that $\zeta, a', d' \in \OF_L^\times$. The upper left $2 \times 2$
block of $A$ is given by
$$
 \mat{\zeta + \alpha z \zeta \varpi^{m}}{z\zeta}{\alpha a'\varpi^{m} -\alpha
 \varpi^{m}(\zeta + \alpha z \zeta \varpi^{m})}
 {\:a'-\alpha z\zeta \varpi^{m}}.$$
We will repeatedly use the following fact:
\begin{equation}\label{well-defined-sub-lemmaeq1}
 \text{If }x \in \OF + \P^r,\text{ then }x \equiv \bar{x} \pmod{(\alpha-\bar\alpha)\P^r}.
\end{equation}
For if $x=y+\alpha z$ with $y\in\OF$ and $z\in\p^r$, then $x-\bar x=(\alpha-\bar\alpha)z$.
Applying this to the matrix entries of $A$, we get $z\zeta \equiv \bar{z} \bar{\zeta} \pmod{(\alpha-\bar\alpha)\P^r}$, and then
\begin{equation}\label{a'-zeta-condition}
 a'-\bar{a'} \equiv (\alpha - \bar{\alpha}) z \zeta \varpi^{m}
 \pmod{(\alpha-\bar\alpha)\P^r}, \qquad \zeta - \bar{\zeta} \equiv
 (\bar{\alpha}-\alpha) z \zeta \varpi^{m} \pmod{(\alpha-\bar\alpha)\P^r}.
\end{equation}
Using $\zeta + \alpha z \zeta \varpi^{m}\equiv \bar\zeta + \bar\alpha \bar{z} \bar\zeta
\varpi^{m} \pmod{(\alpha-\bar\alpha)\P^r}$ and (\ref{a'-zeta-condition}),
we get from the $(2,1)$ coefficient of $A$ that
$$
 (a' \varpi^{m} - \bar\zeta \varpi^{m}) (\alpha - \bar\alpha) \equiv 0
  \pmod{(\alpha-\bar\alpha)\P^r}.
$$
Hence $a' \varpi^{m} - \bar\zeta \varpi^{m}\equiv0\pmod{\P^r}$,
so that $a' \bar{\zeta}^{-1} \in 1+\P^{r-m}$, as required.\qed
\subsubsection*{New- and oldforms in $I(s, \chi, \chi_0, \tau)$}
Let $(\tau,V_\tau)$ be a generic, irreducible, admissible representation of $\GL_2(F)$.
We assume that $V_\tau$ is the Whittaker model of $\tau$ with
respect to the character of $F$ given by $\psi^{-c}(x)=\psi(-cx)$.
Let $\p^n$ be the conductor of $\tau$, where $n$ is a non-negative integer. Let
$W^{(0)}\in V_\tau(n)$ be the local newform as in Lemma \ref{GL2-newform-Kirillov}.
Observe that the central character $\omega_\tau$ is
trivial on $1+\p^n$. We choose any character $\chi_0$ of $L^\times$ such that
\begin{equation}\label{chi0defpropertieseq}
 \chi_0|_{F^\times} = \omega_\tau \qquad \mbox{ and } \qquad\chi_0|_{1+\P^n} = 1.
\end{equation}
(for $n=0$ we mean that $\chi_0$ is unramified).
Given an unramified character $\Lambda$ of $L^\times$, we
define the character $\chi$ of $L^\times$ by the formula
\begin{equation}\label{chidefeq}
 \chi(\zeta) = \Lambda(\bar{\zeta})^{-1} \chi_0(\bar{\zeta})^{-1}.
\end{equation}
Let $I(s, \chi, \chi_0, \tau)$ be the parabolically induced representation of $G(F)$
as defined in Sect.\ \ref{parabolicinductionsec}.
Explicitly, the space of $I(s,\chi,\chi_0,\tau)$ consists of functions
$W:\:G(F)\rightarrow\C$ with the transformation properties
(\ref{Wsharpproperty1eq}) and (\ref{Wsharpproperty2eq}).
The following result shows that there is an essentially unique vector in
$I(s,\chi,\chi_0,\tau)$ right invariant under $K^H\Gamma(\P^n)$.
This vector will be our choice of local section which will be used to evaluate the
non-archimedean local zeta integrals (\ref{localZseq}).
\begin{theorem}\label{unique-W-theorem}
 Let $\chi, \chi_0$ and $\tau$ be as above. Let $\p^n$, $n\geq0$,
 be the conductor of $\tau$. Let
 $$
  V(r) := \{ W \in I(s,\chi,\chi_0,\tau) : W(g \gamma, s) = W(g, s)
   \mbox{ for all } g \in G(F),\:\gamma \in K^H\Gamma(\P^r) \}.
 $$
 Then
 $$
  \dim(V(r)) = \left\{\begin{array}{ll}
   \displaystyle\frac{(r-n+1)(r-n+2)}2& \hbox{ if } r \geq n,\\
   0&\hbox{ if } r < n.\end{array}\right.
 $$
\end{theorem}
{\bf Proof.} Let $W \in V(r)$. By Proposition \ref{disj-doub-coset-decomp-general-n-prop},
$W$ is completely determined by its values on $\eta_m$, $0 \leq m \leq r$.
Let $r\geq m\geq0$.
For any $\mat{a'}{b'}{c'}{d'} \in K^{(0)}(\p^{r-m})$
(see (\ref{K'defeq})), we have
$$
 \begin{bmatrix}1&&&\\&a'&&b'\\&&\mu&\\&c'&&d'\end{bmatrix} \in M(F) N(F) \cap \eta_m K^H
 \Gamma(\P^r) \eta_m^{-1},\quad\mu = a'd'-b'c'.
$$
It follows that
$$
 W(\eta_m) = W(\begin{bmatrix}1&&&\\&a'&&b'\\&&\mu&\\&c'&&d'\end{bmatrix} \eta_m)
  = \tau(\mat{a'}{b'}{c'}{d'})W(\eta_m).
$$
Hence, for $0 \leq m \leq r$, a necessary condition for $v_m:=W(\eta_m)$ is that it is
invariant under $K^{(0)}(\p^{r-m})$.
Since the conductor of $\tau$ is $\p^n$, we conclude that $v_m = 0$ if $r-m < n$.
Therefore $\dim(V(r)) = 0$ for all $r < n$.

Now suppose that $r \geq n$. We will show that, for any $m$ such that $r-m \geq n$,
if $v_m$ is chosen to be any vector in $V_\tau(r-m)$, then we obtain a well-defined
function $W$ in $V(r)$. For $m=r$ this is easy to check, since in this case $n=0$
and all the data is unramified. Assume therefore that $r>m$.
We have to show that for $m_1n_1\eta_m k_1\gamma_1=m_2n_2\eta_m k_2\gamma_2$,
with $m_i\in M(F)$, $n_i\in N(F)$, $k_i\in K^H$ and $\gamma_i\in\Gamma(\P^{r})$,
\begin{equation}\label{Wsharpwelldeflemmanaeq1}
 |N(\zeta_1)\cdot\mu_1^{-1}|^{3(s+1/2)}\chi(\zeta_1)\,
   (\chi_0 \times \tau) (\mat{a'_1}{b'_1}{c'_1}{d'_1}) v_m
 =|N(\zeta_2)\cdot\mu_2^{-1}|^{3(s+1/2)}\chi(\zeta_2)\,
   (\chi_0 \times \tau)(\mat{a'_2}{b'_2}{c'_2}{d'_2}) v_m.
\end{equation}
We have $\eta_m^{-1}m_2^{-1}m_1 n^\ast \eta_m \in K^H \Gamma(\P^{r})$, where
$n^\ast \in N(F)$ depends on $m_1, m_2, n_1, n_2$. Let
$$
 \hat{m}:= m_2^{-1}m_1=\begin{bmatrix}\zeta\\&\tilde a&&\tilde b\\
  &&\mu\bar{\zeta}^{-1}\\&\tilde c&&\tilde d\end{bmatrix},
 \quad n^\ast = \begin{bmatrix}
                 1 & z &  &  \\
                  & 1 &  &  \\
                  &  & 1 &  \\
                  &  & -\overline{z} & 1 \\
               \end{bmatrix}
               \begin{bmatrix}
                 1 &  & w & y \\
                  & 1 & \overline{y} &  \\
                  &  & 1 &  \\
                  &  &  & 1 \\
               \end{bmatrix}.
$$
Then $\zeta\in\OF_L^\times$ and $\mu\in\OF^\times$.
By definition, $\zeta_1=\zeta_2\zeta$ and $\mu_1=\mu_2\mu$. Hence
(\ref{Wsharpwelldeflemmanaeq1}) is equivalent to
\begin{equation}\label{Wsharpwelldeflemmanaeq2}
 \chi(\zeta)\,(\chi_0 \times \tau)(\mat{a'_1}{b'_1}{c'_1}{d'_1}) v_m
 = (\chi_0 \times \tau)(\mat{a'_2}{b'_2}{c'_2}{d'_2}) v_m.
\end{equation}
Using Lemma \ref{well-defined-sub-lemma}, we get $\tilde a\bar\zeta^{-1}\in 1+\P^{r-m}$ and
$\tilde c\in \P^{r-m}$. Hence, using (\ref{chi0defpropertieseq}) and (\ref{chidefeq})
(with unramified $\Lambda$) and the fact that $v_m \in V_\tau(r-m)$,
\begin{align*}
 \chi(\zeta)\,(\chi_0 \times \tau)(\mat{a'_1}{b'_1}{c'_1}{d'_1})v_m
  &=\chi(\zeta)\,(\chi_0 \times \tau) (\mat{a'_2}{b'_2}{c'_2}{d'_2}
   \mat{\tilde a}{\tilde b}{\tilde c}{\tilde d}) v_m\\
  &=\chi(\zeta)\chi_0(\tilde a)\,(\chi_0 \times \tau)(\mat{a'_2}{b'_2}{c'_2}{d'_2}
   \mat{1}{\tilde b/\tilde a}{\tilde c/\tilde a}{\tilde d/\tilde a})v_m\\
  &=\chi_0(\bar\zeta^{-1})\chi_0(\tilde a)\,(\chi_0 \times \tau)(\mat{a'_2}{b'_2}{c'_2}{d'_2})v_m\\
  &=(\chi_0 \times \tau)(\mat{a'_2}{b'_2}{c'_2}{d'_2})v_m,
\end{align*}
as claimed.

Now, using the formula for $\dim(V_\tau(r-m))$ from Theorem \ref{GL2newformtheorem}
completes the proof of the theorem. \qed
\subsection{The zeta integral}\label{padiczetasec}
As in the previous section let $(\tau,V_\tau)$ be a generic, irreducible, admissible
representation of $\GL_2(F)$ with conductor $\p^n$.
We assume that $V_\tau$ is the Whittaker model of $\tau$ with
respect to the additive character $\psi^{-c}(x)=\psi(-cx)$. Let the characters
$\chi_0$ and $\chi$ of $L^\times$ be as in (\ref{chi0defpropertieseq}),
resp.\ (\ref{chidefeq}). In the induced representation $I(s,\chi,\chi_0,\tau)$,
consider the spaces $V(r)$ of invariant vectors defined in Theorem \ref{unique-W-theorem}.
Taking $r=n$ in this theorem, we see that $\dim(V(n)) = 1$. The proof of
Theorem \ref{unique-W-theorem} shows that, in the model
$I_W(s,\chi,\chi_0,\tau)$ of $I(s,\chi,\chi_0,\tau)$ consisting of
complex-valued functions (see Sect.\ \ref{parabolicinductionsec}),
$V(n)$ is spanned by the unique function $W^\#(\,\cdot\,,s)$ with the following properties.
 \begin{itemize}
  \item If $g\notin M(F)N(F)\eta K^H \Gamma(\P^{n})$, then $W^\#(g,s)=0$.
  \item If $g=mn\eta k \gamma$ with $m\in M(F)$, $n\in N(F)$, $k\in K^H$,
   $\gamma \in \Gamma(\P^{n})$, then $W^\#(g,s)=W^\#(m \eta,s)$.
  \item For $\zeta\in L^\times$ and $\mat{a'}{b'}{c'}{d'}\in\GU(1,1;L)(F)$,
   \begin{equation}\label{Wsharpformulaeq}
     W^\#(\begin{bmatrix}\zeta\\&1\\&&\bar{\zeta}^{-1}\\&&&1\end{bmatrix}
    \begin{bmatrix}1\\&a'&&b'\\&&\mu\\&c'&&d'\end{bmatrix} \eta,s)
    =|N(\zeta)\cdot\mu^{-1}|^{3(s+1/2)}\chi(\zeta)\,
    W^{(0)}(\mat{a'}{b'}{c'}{d'}).
   \end{equation}
   Here $\mu=\bar{a'}d'-b'\bar{c'}$ and $W^{(0)}$ is the newform in $\tau$
   as defined in Lemma \ref{GL2-newform-Kirillov}, but extended to a function
   on $\GU(1,1;L)(F)$ via the character $\chi_0$ as in (\ref{extendedWformulaeq}).
\end{itemize}
It is this function $W^\#$ for which we will evaluate the local integral
$Z(s,W^\#,B)$ defined in (\ref{localZseq}). The other ingredient in this integral
is the Bessel function $B$, which is the spherical vector in the Bessel model
of an unramified representation $(\pi,V_\pi)$ of $H(F)$
with respect to the character $\Lambda\otimes\theta$ of $R(F)$; see Sect.\ \ref{sphericalbesselfunctionsec}. In the following we shall assume $n>0$, since
for unramified $\tau$ the local integral has been computed by Furusawa; see
Theorem (3.7) in \cite{Fu}. Since both functions $B$ and $W^\#$ are right
$K^H$-invariant, it follows from (\ref{RFKHrepresentativeseq}) that
the integral (\ref{localZseq}) is given by
\begin{equation}\label{zeta-integral-sum}
 Z(s,W^\#,B)=\sum\limits_{l,m \geq 0} B(h(l,m))W^\#(\eta h(l,m),s)V_mq^{3m+3l}.
\end{equation}
Here, as in Sect.\ 3.5 of \cite{Fu},
$$
 V_m = \int\limits_{T(F)\backslash T(F)
 \left[\begin{smallmatrix}\varpi^m&\\&1\end{smallmatrix}\right] \GL_2(\OF)} dt.
$$
We will only need the value of $V_0$, which is normalized to be equal to $1$.
To compute the integral (\ref{zeta-integral-sum}), we need to know for what values
of $l,m$ does $\eta h(l,m)$ belong to the support of $W^\#$.
Since $\eta h(l,m) = h(l,m) \eta_m$, with $h(l,m) \in M(F)$ and
$\eta_m$ as in (\ref{etamdef2eq}), all that is relevant is
for what values of $m$ is $\eta_m$ in the support of $W^\#$.
The support of $W^\#$ is $P(F)\eta K^H\Gamma(\P^n)=P(F)\eta K^H$
(see Lemma \ref{KHGammaKHlemma}). Hence, by Proposition
\ref{disj-doub-coset-decomp-general-n-prop} iii), only $\eta_0=\eta$ is
in the support. It follows that the integral (\ref{zeta-integral-sum}) reduces to
\begin{equation}\label{zeta-integral-sum-m=0}
 Z(s,W^\#,B)= \sum\limits_{l \geq 0} B(h(l,0))W^\#(\eta h(l,0),s)q^{3l}.
\end{equation}
By  (\ref{chidefeq}) and (\ref{Wsharpformulaeq}),
\begin{align}\label{W-value-h(l,0)}
 W^\#(\eta h(l,0),s)&= |N(\varpi^l)\varpi^{-l}|^{3(s+1/2)}
  \chi(\varpi^l) W^{(0)}(\mat{\varpi^l}{}{}{1})\nonumber\\
 &=q^{-3(s+1/2)l} \omega_\pi(\varpi^{-l}) \omega_\tau(\varpi^{-l})
  W^{(0)}(\mat{\varpi^l}{}{}{1}).
\end{align}
We will consider three cases for the representation $\tau$ according to the values of the newform $W^{(0)}$ given in Lemma \ref{GL2-newform-Kirillov}.

{\bf Case 1:} Let $\tau$ be either a supercuspidal representation,
or a twist $\Omega\St_{\GL(2)}$ of the Steinberg representation with
a ramified character $\Omega$, or a principal series representation $\alpha\times\beta$
with two ramified characters $\alpha, \beta$ (such that $\alpha\beta^{-1}\neq|\,|^{\pm1}$).
In each of these cases, using Lemma \ref{GL2-newform-Kirillov} i), we have
\begin{equation}\label{zeta-integral-final-1}
 Z(s,W^\#,B)=1.
\end{equation}

{\bf Case 2:} Let $\tau=\alpha\times\beta$ be a principal series representation
with an unramified character $\alpha$ and a ramified character $\beta$.
Then, by Lemma \ref{GL2-newform-Kirillov} ii) and (\ref{W-value-h(l,0)}),
\begin{align*}
 Z(s,W^\#,B)&= \sum\limits_{l \geq 0} B(h(l,0)) q^{-3(s+\frac 12)l} \omega_\pi(\varpi^{-l})
  \omega_\tau(\varpi^{-l})(\beta(\varpi)q^{-1/2})^l q^{3l} \\
 &=\sum\limits_{l \geq 0}B(h(l,0)) \Big(q^{-3s+1}(\omega_\pi \alpha)^{-1}(\varpi)\Big)^l.
\end{align*}
Let $\varpi_L$ be the uniformizer of $\OF_L$ if $L$ is a field, and set
$\varpi_L = (\varpi,1)$ if $L$ is not a field. If $L/F$ is a ramified field
extension, we assume in addition that $N_{L/F}(\varpi_L)=\varpi$.
Then, using the notations from Sect.\ \ref{sphericalbesselfunctionsec},
\begin{equation}\label{zeta-integral-final-2}
 Z(s,W^\#,B)=\frac{H(y)}{Q(y)}\qquad\text{with}\qquad
 y=q^{-3s+1}(\omega_\pi \alpha)^{-1}(\varpi).
\end{equation}
Explicitly,
\begin{align}\label{Q-formula}
 Q(y)&=\prod_{i=1}^4\big(1-\gamma^{(i)}(\varpi)q^{-3/2}q^{-3s+1}
   (\omega_\pi\alpha)^{-1}(\varpi)\big) \nonumber\\
 &=\prod_{i=1}^4\big(1-q^{-3s-\frac 12}(\gamma^{(i)}\alpha)^{-1}(\varpi)\big),
\end{align}
and
\begin{equation}\label{H-formula}
 H(y) = \left\{\begin{array}{ll}
    1-\big(\Lambda (\omega_\pi \alpha)^{-2}\big)(\varpi)q^{-6s-2}& \hbox{ if }
     \big(\frac L\p\big)=-1,\\[1ex]
    1-\Lambda(\varpi_L) (\omega_\pi \alpha)^{-1}(\varpi)q^{-3s-1}& \hbox{ if }
     \big(\frac L\p\big)= 0,\\[1ex]
    (1-\Lambda(\varpi_L) (\omega_\pi \alpha)^{-1}(\varpi)q^{-3s-1})
    (1-\Lambda(\varpi\varpi_L^{-1})(\omega_\pi \alpha)^{-1}(\varpi)q^{-3s-1})
     & \hbox{ if } \big(\frac L\p\big)= 1.
 \end{array}\right.
\end{equation}

{\bf Case 3:} Let $\tau=\Omega\St_{\GL(2)}$ with an unramified character $\Omega$
of $F^\times$. Then, using Lemma \ref{GL2-newform-Kirillov} iii) and (\ref{W-value-h(l,0)}),
a similar calculation as in Case 2 shows that
\begin{equation}\label{zeta-integral-final-3}
 Z(s,W^\#,B)=\frac{H(y)}{Q(y)}\qquad\text{with}\qquad
 y=q^{-3s+1/2}(\omega_\pi\Omega)^{-1}(\varpi).
\end{equation}

We now get the following theorem, which is our main non-archimedean result.
\begin{theorem}\label{ram-central-char-main-thm-local}
 Let $\pi$ be an irreducible, admissible, unramified representation
 of $\GSp_4(F)$ and let $\tau$ be an irreducible, admissible
 representation of $\GL_2(F)$. Let $B$ be the unramified Bessel function given
 by formula (\ref{suganox0eq}). Let $W^\#$ be the element of $I_W(s,\chi,\chi_0,\tau)$
 defined in Sect.\ \ref{padiczetasec}. Then the local zeta
 integral $Z(s,W^\#,B)$ defined in (\ref{localZseq}) is given by
 \begin{equation}\label{final-integral-l-fn-formula}
  Z(s,W^\#,B)=\frac{L(3s+\frac 12, \tilde\pi \times \tilde\tau)}
  {L(6s+1,\chi|_{F^\times})L(3s+1,\tau \times \AI(\Lambda) \times \chi|_{F^\times})}Y(s),
 \end{equation}
 where
 \begin{align*}
  Y(s)=\left\{\begin{array}{l@{\qquad}l}
   1&\text{if }\tau=\alpha\times\beta,\;\alpha,\beta\text{ unramified},\\
   L(6s+1,\chi|_{F^\times})
    &\text{if }\tau=\alpha\times\beta,\;\alpha\text{ unram.},\:\beta\text{ ram.},\:
     \big(\frac L\p\big)=\pm1,\\
&\qquad\text{OR }\tau=\alpha\times\beta,\;\alpha\text{ unram.},\:\beta\text{ ram.},\\
      &\qquad\qquad\big(\frac L\p\big)=0 \text{ and }\beta\chi_{L/F}\text{ ramified},\\
&\qquad\text{OR } \tau=\Omega\St_{\GL(2)},\;\Omega\text{ unramified},\\
     \displaystyle\frac{L(6s+1,\chi|_{F^\times})}
      {1-\Lambda(\varpi_L)(\omega_\pi \beta)^{-1}(\varpi)q^{-3s-1}}
    &\text{if }\tau=\alpha\times\beta,\;\alpha\text{ unram.},\:\beta\text{ ram.},\:
     \big(\frac L\p\big)=0,\\
    &\qquad\text{and }\beta\chi_{L/F}\text{ unramified},\\
   L(6s+1,\chi|_{F^\times})L(3s+1,\tau \times \AI(\Lambda) \times \chi|_{F^\times})
    &\text{if }\tau=\alpha\times\beta,\;\alpha,\beta\text{ ramified},\\
    &\qquad\text{OR } \tau=\Omega\St_{\GL(2)},\:\Omega\text{ ramified},\\
    &\qquad\text{OR } \tau \text{ supercuspidal}.
  \end{array}\right.
 \end{align*}
 In (\ref{final-integral-l-fn-formula}), $\tilde\pi$ and $\tilde\tau$ denote the
 contragredient of $\pi$ and $\tau$, respectively. The symbol $\AI(\Lambda)$
 stands for the $\GL_2(F)$ representation attached to the character $\Lambda$
 of $L^\times$ via automorphic induction, and
 $L(3s+1,\tau \times \AI(\Lambda) \times \chi|_{F^\times})$ is
 a standard $L$-factor for $\GL_2\times\GL_2\times\GL_1$.
\end{theorem}
{\bf Proof.} If $\tau=\alpha\times\beta$ with unramified $\alpha$ and $\beta$, then
this is Theorem (3.7) in Furusawa's paper \cite{Fu}. If $\tau=\alpha\times\beta$
with unramified $\alpha$ and ramified $\beta$ (Case 2 above), then, from
the local Langlands correspondence, we have the following
$L$-functions attached to the representations $\tilde\pi \times
\tilde\tau$ of $\GSp_4(F) \times \GL_2(F)$ and $\tau \times
\AI(\Lambda) \times \chi|_{F^\times}$ of $\GL_2(F) \times \GL_2(F)\times \GL_1(F)$,
\begin{equation}\label{gsp4-gl2-l-fn}
 L(s, \tilde\pi \times \tilde\tau)=\prod_{i=1}^4
 \big(1-(\gamma^{(i)}\alpha)^{-1}(\varpi)q^{-s}\big)^{-1}
\end{equation}
and
\begin{equation}\label{gl2-gl2-gl1-l-fn}
 \frac1{L(s, \tau \times \AI(\Lambda) \times \chi|_{F^\times})}=
  \left\{\begin{array}{l@{\;}l}
    1-\big(\Lambda (\omega_\pi \alpha)^{-2}\big)(\varpi)q^{-2s}
     &\text{if } \big(\frac L\p\big)=-1,\\
    1-\Lambda(\varpi_L) (\omega_\pi \alpha)^{-1}(\varpi)q^{-s}
     &\text{if } \big(\frac L\p\big)= 0\text{ and}\\
     &\beta\chi_{L/F}\text{ ram.},\\
    (1-\Lambda(\varpi_L) (\omega_\pi \alpha)^{-1}(\varpi)q^{-s})
     (1-\Lambda(\varpi_L) (\omega_\pi \beta)^{-1}(\varpi)q^{-s})
     &\text{if } \big(\frac L\p\big)=0\text{ and}\\
     &\beta\chi_{L/F}\text{ unram.},\\
    (1-\Lambda(\varpi_L) (\omega_\pi \alpha)^{-1}(\varpi)q^{-s})
     (1-\Lambda(\varpi\varpi_L^{-1})(\omega_\pi \alpha)^{-1}(\varpi)q^{-s})
     &\text{if } \big(\frac L\p\big)= 1.\end{array}\right.
\end{equation}
The desired result therefore follows from (\ref{Q-formula}) and (\ref{H-formula}).
If $\tau$ is an unramified twist of the Steinberg representation
(Case 3 above), then the result was proved in Theorem 3.8.1 of \cite{PS1}.
In all remaining cases (i.e., Case 1 above) we have $L(s,\tilde\pi\times\tilde\tau)=1$,
so that the theorem follows from (\ref{zeta-integral-final-1}).
This completes the proof.\qed
\section{Local archimedean theory}
In this section we evaluate the local zeta integral (\ref{localZseq}) in the
real case. As in the non-archimedean case, the key steps are the choices of the
vector $W^\#$ and the actual computation of the integral $Z(s,W^\#,B)$. 
\subsection{Notations}\label{arch-notations}
We recall some of the definitions and basic facts from Sect.\ 4.1 of \cite{PS1}.
Let $G=\GU(2,2;\C)$ as in Sect.\ \ref{unitarygroupsec} (with $F=\R$ and $L=\C$).
Consider the symmetric domains $\HH_2 := \{Z \in M_2(\C):\;
i(\,^{t}\!\bar{Z}-Z) \mbox{ is positive definite}\}$ and $\SH_2 :=
\{Z \in \HH_2:\;^{t}Z = Z \}$. The group $G^{+}(\R) := \{g \in
G(\R) : \mu_2(g) > 0 \}$ acts on $\HH_2$ via $(g,Z) \mapsto
g\langle Z \rangle$, where
$$
 g\langle Z \rangle = (AZ+B)(CZ+D)^{-1}, \mbox{ for } g =
 \mat{A}{B}{C}{D} \in G^{+}(\R), Z \in \HH_2.
$$
Under this action, $\SH_2$ is stable by $H^+(\R) = \GSp_4^+(\R)$.
The group $K_\infty=\{g\in G^{+}(\R): \mu_2(g) = 1,\, g\langle
I\rangle = I \}$ is a maximal compact subgroup of $G^{+}(\R)$.
Here, $I=\mat{i}{}{}{i} \in\HH_2$. Explicitly,
$$
 K_\infty=\{\mat{A}{B}{-B}{A}:\;A,B\in M(2,\C),\;^t\!\bar AB=\,^t\bar BA,\;
 ^t\!\bar AA+\,^t\bar BB=1\}.
$$
By the Iwasawa decomposition
\begin{equation}\label{realiwasawaeq}
 G(\R)=M^{(1)}(\R)M^{(2)}(\R)N(\R)K_\infty,
\end{equation}
where $M^{(1)}(\R)$, $M^{(2)}(\R)$ and $N(\R)$ are as defined in
(\ref{M1defn}), (\ref{M2defn}) and (\ref{Ndefn}). A calculation shows that
\begin{equation}\label{realKintersect1eq}
 M^{(1)}(\R)M^{(2)}(\R)N(\R)\cap K_\infty
 =\{\begin{bmatrix}\zeta\\&\alpha&&\beta\\&&\zeta\\&-\beta&&\alpha\end{bmatrix}:
 \:\zeta,\alpha,\beta\in\C,\:|\zeta|=1,\:|\alpha|^2+|\beta|^2=1,\:
 \alpha\bar\beta=\beta\bar\alpha\}.
\end{equation}
Note also that
\begin{equation}\label{realKintersect2eq}
 M^{(2)}(\R)\cap K_\infty
 =\{\begin{bmatrix}1\\&\alpha&&\beta\\&&1\\&-\beta&&\alpha\end{bmatrix}:
 \:\alpha,\beta\in\C,\:|\alpha|^2+|\beta|^2=1,\:\alpha\bar\beta=\beta\bar\alpha\},
 \end{equation}
and that there is an isomorphism
\begin{align}\label{realKintersect3eq}
 (S^1\times\SO(2))/\{(\lambda,\mat{\lambda}{}{}{\lambda}):\:\lambda=\pm1\}&\stackrel{\sim}{\longrightarrow}
  M^{(2)}(\R)\cap K_\infty,\nonumber\\
 (\lambda,\mat{\alpha}{\beta}{-\beta\:}{\alpha})&\longmapsto
  \begin{bmatrix}1\\&\lambda\alpha&&\lambda\beta\\&&1\\&-\lambda\beta&&\lambda\alpha\end{bmatrix}.
\end{align}
For $g\in G^+(\R)$ and $Z\in\mathbb{H}_2$, let $J(g,Z)=CZ+D$ be
the automorphy factor. Then, for any integer $l$, the map
\begin{equation}\label{realcompactcharactereq}
 k\longmapsto \det(J(k,I))^l
\end{equation}
defines a character $K_\infty\rightarrow\C^\times$. If $k\in
M^{(2)}(\R)\cap K_\infty$ is written in the form
(\ref{realKintersect3eq}), then
$\det(J(k,I))^l=\lambda^le^{-il\theta}$, where
$\alpha=\cos(\theta)$, $\beta=\sin(\theta)$. Let $K^H_\infty =
K_{\infty} \cap H^+(\R)$. Then $K^H_\infty$ is a maximal compact
subgroup, explicitly given by
$$
 K^H_\infty=\{\mat{A}{B}{-B}{A}:\;^t\!AB=\,^tBA,\;^t\!AA+\,^tBB=1\}.
$$
Sending $\mat{A}{B}{-B}{A}$ to $A-iB$ gives an isomorphism
$K^H_\infty\cong{\rm U}(2)$. Recall that we have chosen
$a,b,c\in\R$ such that $d=b^2-4ac\neq0$. In the archimedean case
we shall assume that $d<0$ and let $D=-d$. Then
$\R(\sqrt{-D})=\C$. The group $T(\R)$ defined in (\ref{TFdefeq}) is given by
\begin{equation}\label{TReq}
 T(\R)=\{\mat{x+yb/2}{yc}{-ya}{x-yb/2}:\:x,y\in\R,\:x^2+y^2D/4>0\}.
\end{equation}
Let
\begin{equation}\label{Tinfty1eq}
 T^1(\R)=T(\R)\cap\SL(2,\R)
 =\{\mat{x+yb/2}{yc}{-ya}{x-yb/2}:\:x,y\in\R,\:x^2+y^2D/4=1\}.
\end{equation}
We have $T(\R)\cong\C^\times$ via
$\mat{x+yb/2}{yc}{-ya}{x-yb/2}\mapsto x+y\sqrt{-D}/2$. Under this
isomorphism $T^1(\R)$ corresponds to the unit circle. We have
\begin{equation}\label{TRTinfty1eq}
 T(\R)=T^1(\R)\cdot\{\mat{\zeta}{}{}{\zeta}:\:\zeta>0\}.
\end{equation}
As in \cite{Fu}, p.\ 211, let $t_0\in\GL_2(\R)^+$ be such that
$T^1(\R)=t_0\SO(2)t_0^{-1}$. We will make a specific choice of
$t_0$ when we choose the matrix $S=\mat{a}{b/2}{b/2}{c}$ below. It
is not hard to see that
\begin{equation}\label{HRBesseldecompeq}
 H(\R)=R(\R)\cdot\big\{\begin{bmatrix}\lambda t_0\mat{\zeta}{}{}{\zeta^{-1}}&\\
  &^tt_0^{-1}\mat{\zeta^{-1}}{}{}{\zeta}\end{bmatrix}:\:\lambda\in\R^\times,\,\zeta\geq1\big\}
  \cdot K^H_\infty.
\end{equation}
Here, $R(\R)=T(\R)U(\R)$ is the Bessel subgroup. One can check
that all the double cosets in (\ref{HRBesseldecompeq}) are
disjoint.

\subsection{The Bessel function}\label{realbesselfnsec}
Recall that we have chosen three elements $a,b,c\in\R$ such that
$d=b^2-4ac<0$. We will now make the stronger assumption that
$S=\mat{a}{b/2}{b/2}{c}\in M_2(\R)$ is a positive definite matrix.
Set $D=4ac-b^2>0$, as above. Given a positive integer $l\geq2$,
consider the function $B:\:H(\R)\rightarrow\C$ defined by
\begin{equation}\label{archBesselformula2eq}\renewcommand{\arraystretch}{1.2}
 B(h) := \left\{\begin{array}{ll}
 \mu_2(h)^l\,\overline{\det(J(h, I))^{-l}}\,e^{-2\pi i\,
 {\rm tr}(S\overline{h\langle I\rangle})}& \hbox{ if } h \in H^+(\R),\\
 0& \hbox{ if } h \notin H^+(\R),\\
\end{array}\right.
\end{equation}
where $I=\mat{i}{}{}{i}$. Note that the function $B$ only depends
on the choice of $S$ and $l$. Recall the character $\theta$ of
$U(\R)$ defined by $\theta(\mat{1}{X}{}{1}) = \psi({\rm tr}(SX))$.
It depends on the choice of additive character $\psi$, and
throughout we choose $\psi(x)=e^{-2\pi ix}$. Then the function $B$
satisfies
\begin{equation}\label{realBproperty1eq}
 B(tuh)=\theta(u)B(h)\qquad\text{for }
 h\in H(\R),\;t\in T(\R),\;u\in U(\R),
\end{equation}
and
\begin{equation}\label{realBproperty2eq}
 B(hk)=\det(J(k,I))^lB(h)\qquad\text{for }
 h\in H(\R),\;k\in K^H_\infty.
\end{equation}
Property (\ref{realBproperty1eq}) means that $B$ satisfies the
Bessel transformation property with the character $\Lambda \otimes
\theta$ of $R(\R)$, where $\Lambda$ is trivial. In fact, by the
considerations in \cite{Su} 1-3, or by \cite{PS} Theorem 3.4, $B$
is the highest weight vector (weight $(-l,-l)$) in a holomorphic discrete series
representation (or limit of such if $l=2$) of ${\rm PGSp}_4(\R)$
corresponding to Siegel modular forms of degree $2$ and weight
$l$. By (\ref{realBproperty1eq}) and (\ref{realBproperty2eq}), the
function $B$ is determined by its values on a set of
representatives for $R(\R)\backslash H(\R)/K^H_\infty$. Such a set
is given in (\ref{HRBesseldecompeq}).
\subsection{The function $W^\#$}\label{realW sharp function}
Let $(\tau,V_\tau)$ be a  generic, irreducible, admissible
representation of $\GL_2(\R)$ with central character
$\omega_{\tau}$. We assume that
$V_\tau=\mathcal{W}(\tau,\psi_{-c})$ is the Whittaker model of
$\tau$ with respect to the non-trivial additive character
$x\mapsto\psi(-cx)$. Note that $S$ positive definite implies
$c>0$. Let $W^{(0)}\in V_\tau$ have weight $l_1$. Then $W^{(0)}$
has the properties
\begin{equation}\label{W0prop1eq}
 W^{(0)}(gr(\theta))=e^{il_1\theta}W^{(0)}(g),\qquad g\in\GL_2(\R),\;r(\theta)
 =\mat{\cos(\theta)}{\sin(\theta)}{-\sin(\theta)}{\cos(\theta)}\in\SO(2),
\end{equation}
and
\begin{equation}\label{W0prop2eq}
 W^{(0)}(\mat{1}{x}{}{1}g)=\psi(-cx)W^{(0)}(g),\qquad g\in\GL_2(\R),\;x\in\R.
\end{equation}
Let $l_2$ be an integer of the same parity as $l_1$; further below
in (\ref{l2defeq}) we will be more
specific. Let $\chi_0$ be the character of $\C^\times$ with the properties
\begin{equation}\label{realchi0defeq}
 \chi_0\big|_{\R^\times}=\omega_\tau,\qquad\chi_0(\zeta)=\zeta^{-l_2}\quad\text{for }
 \zeta\in\C^\times,\:|\zeta|=1.
\end{equation}
Such a character exists since $\omega_\tau(-1)=(-1)^{l_1}=(-1)^{l_2}$. We
extend $W^{(0)}$ to a function on $M^{(2)}(\R)$ via
\begin{equation}\label{realW1extensioneq}
 W^{(0)}(\zeta g)=\chi_0(\zeta)W^{(0)}(g),\qquad \zeta\in\C^\times,\:g\in\GL_2(\R);
\end{equation}
see (\ref{extendedWformulaeq}).
We will need the values of $W^{(0)}$ at elements $\mat{t}{}{}{1}$ for
$t \neq 0$. For this we consider the Lie algebra
$\mathfrak{g}=\mathfrak{gl}(2,\R)$ and its elements
$$
 R=\mat{0}{1}{0}{0},\qquad L=\mat{0}{0}{1}{0},\qquad
 H=\mat{1}{0}{0}{-1},\qquad Z=\mat{1}{0}{0}{1}.
$$
In the universal enveloping algebra $U(\mathfrak{g})$ let $\Delta=\frac14(H^2+2RL+2LR)$.
Then $\Delta$ lies in the center of $U(\mathfrak{g})$ and acts on
$V_\tau$ by a scalar, which we write in the form $-(\frac14+(\frac r2)^2)$ with $r\in\C$.
In particular,
\begin{equation}\label{W1Deltapropertyeq}
 \Delta W^{(0)}=-\Big(\frac14+\Big(\frac r2\Big)^2\Big)W^{(0)}.
\end{equation}
If one restricts the function $W^{(0)}$ to $\mat{t^{1/2}}{}{}{t^{-1/2}}$, $t>0$, then
(\ref{W1Deltapropertyeq}) reduces to the differential equation
satisfied by the classical Whittaker functions. Hence,
there exist constants $a^+,a^-\in\C$ such that
\begin{equation}\label{W1Whittakerformulaeq}
 W^{(0)}(\mat{t}{0}{0}{1}) = \renewcommand{\arraystretch}{1.3}
  \left\{\begin{array}{ll}
 a^+\omega_\tau((4\pi ct)^{1/2})W_{\frac{l_1}2, \frac{ir}2}(4\pi ct) & \hbox{ if } t > 0, \\
 a^-\omega_\tau((-4\pi ct)^{1/2})
  W_{-\frac{l_1}2, \frac{ir}2}(-4\pi ct)& \hbox{ if } t <0 .\end{array}\right.
\end{equation}
Here, $W_{\pm\frac{l_1}2, \frac{ir}2}$ denotes a classical
Whittaker function; see \cite[p.\ 244]{Bu}, \cite{MOS}. According to
(\ref{chi-lambda-char-condition}), we let $\chi$ be the character of $\C^\times$ given by
\begin{equation}\label{archchidefeq}
 \chi(\zeta)=\chi_0(\bar{\zeta})^{-1}.
\end{equation}
We interpret $\chi$ as a character of $M^{(1)}(\R)\cong\C^\times$.
We wish to define a function $W^\#$ of the form
\begin{equation}\label{Wsharpwelldeflemmaeq1}
  W^\#(m_1m_2nk,s)=\delta_P^{s+1/2}(m_1m_2)f(k)\chi(m_1)W^{(0)}(m_2),
\end{equation}
where $m_1\in M^{(1)}(\R)$, $m_2\in M^{(2)}(\R)$, $n\in N(\R)$ and
$k\in K_\infty$, for some analytic function $f$ on $K_\infty$. Any such
$W^\#$ would be a legitimate section of the induced representation
$I_\C(s,\chi,\chi_0,\tau)$ considered in Sect.\ \ref{parabolicinductionsec}.
In addition, we would like $W^\#$ to satisfy the right transformation property
\begin{equation}\label{Wsharprighttransformeq}
 W^\#(gk,s)=\det(J(k,I))^{-l}W^\#(g,s)\qquad\text{for }
 g\in G(\R),\;k\in K_\infty^H.
\end{equation}
We need this property so that the function $B(g)W^\#(g,s)$
will be right invariant under $K^H_\infty$; see (\ref{realBproperty2eq}).
The following lemma gives the precise conditions to be satisfied by the function
$f$ so that $W^\#$ is well-defined.
\begin{lemma}\label{realWsharpwelldeflemma}
 Let $f$ be a function on $K_\infty$. For $\zeta \in \C$, set $\hat\zeta_1 = \begin{bmatrix}\zeta\\&1\\&&\zeta\\&&&1\end{bmatrix}$ and $\hat\zeta_2 = \begin{bmatrix}1\\&\zeta\\&&1\\&&&\zeta\end{bmatrix}$. For $\theta \in \R$, set $\hat{r}(\theta) = \begin{bmatrix}1\\&\cos(\theta)&&\sin(\theta)\\&&1\\
    &-\sin(\theta)&&\cos(\theta)\end{bmatrix}$.
We can define a function $W^\#(\,\cdot\,,s)$ on $G(\R)$ by formula
(\ref{Wsharpwelldeflemmaeq1}) if and only if $f$ satisfies, for all
$k\in K_\infty$, $\zeta\in S^1$ and $\theta \in \R$, the conditions
   \begin{equation}\label{Wsharpwelldeflemmaeq2}
    f(\hat\zeta_1k)=\chi(\zeta)f(k), \qquad
    f(\hat\zeta_2k)=\chi_0(\zeta)f(k),
   \end{equation}
  \begin{equation}\label{Wsharpwelldeflemmaeq4}
    f(\hat{r}(\theta)k)=e^{il_1\theta}f(k).
   \end{equation}
\end{lemma}
{\bf Proof.} This is obtained by direct computation.\qed

We will now demonstrate how to obtain a function $f$ on $K_\infty$ satisfying all the
required conditions. We define four
functions $\hat a$, $\hat b$, $\hat c$, $\hat d$ on $K_\infty$ by
\begin{align*}
& \hat a(g)=\text{$a$-coefficient of }J(g\,^t\!g,I), \qquad
 \hat b(g)=\text{$b$-coefficient of }J(g\,^t\!g,I),\\
& \hat c(g)=\text{$c$-coefficient of }J(g\,^t\!g,I), \qquad
 \hat d(g)=\text{$d$-coefficient of }J(g\,^t\!g,I).
\end{align*}
The function $J(g,Z)$ was defined in Sect.\ \ref{arch-notations}.
Here, we have written $J(g\,^t\!g,I)$ as $\mat{a}{b}{c}{d}$. Since $h\,^th=1$ for all
$h\in K^H_\infty$, each of these functions is right $K^H_\infty$
invariant.  Calculations show that
\begin{align*}
& \hat a(\hat\zeta_1g)=\zeta^2\hat a(g), \quad
  \hat a(\hat\zeta_2g)=\hat a(g), \quad
 \hat b(\hat\zeta_1g)=\zeta\hat b(g), \quad
  \hat b(\hat\zeta_2g)=\zeta\hat b(g),\\
& \hat c(\hat\zeta_1g)=\zeta\hat c(g),\quad
  \hat c(\hat\zeta_2g)=\zeta\hat c(g),\quad
 \hat d(\hat\zeta_1g)=\hat d(g),\quad
  \hat d(\hat\zeta_2g)=\zeta^2\hat d(g).
\end{align*}
for $\zeta\in S^1$, as well as
$$
 \hat a(\hat{r}(\theta)g)=\hat a(g), \quad
\hat b(\hat{r}(\theta)g)=e^{i\theta}\hat b(g), \quad \hat c(\hat{r}(\theta)g)=e^{-i\theta}\hat c(g), \quad
\hat d(\hat{r}(\theta)g)=\hat d(g)
$$
for all $\theta\in\R$.
\begin{lemma}\label{fonKinftylemma}
 Let $l$ be any integer. Let $t_1,t_2,t_3$ be integers of the same parity satisfying
 $t_1\geq t_2\geq-t_3\geq-2l-t_2$. Then
 there exists a real analytic function $f$ on $K_\infty$ with the following properties.
 \begin{itemize}
  \item For all $h\in K^H_\infty$,
   $$
    f(gh)=\det(J(h,I))^{-l}f(g).
   $$
  \item For all $\zeta\in S^1$,
   $$
    f(\hat\zeta_1g)=\zeta^{t_1}f(g), \qquad
     f(\hat\zeta_2g)=\zeta^{t_2}f(g).
   $$
  \item For all $\theta\in\R$,
   $$
    f(\hat{r}(\theta)g)=e^{it_3\theta}f(g).
   $$
 \end{itemize}
 In fact, for any integer $t\geq0$ such that all exponents in the function
 $$
  f(g)=\hat a(g)^{\frac{t_1-t_2}2+t}\,\hat b(g)^{\frac{t_2+t_3}2-t}
   \,\hat c(g)^{\frac{t_2-t_3}2+l-t}\,\hat d(g)^t\det(J(g,I))^{-l}
 $$
 are non-negative, this function has the desired properties.
\end{lemma}
{\bf Proof.} This follows from the above transformation properties
of $\hat a, \hat b, \hat c, \hat d$. In order to obtain a well-defined function,
we need to make sure all exponents are non-negative integers.\qed

With $l$ being the weight of our Siegel modular form and $l_1$ being the
weight of our function $W^{(0)}\in V_\tau$, we now make the choice $t=0$ and
$t_3=l_1$ in the above lemma, and $t_1=t_2$ as small as possible. We obtain the analytic
function on $K_\infty$ given by
\begin{equation}\label{finalfdefeq}
 f(g)=\left\{\begin{array}{l@{\qquad\text{if }}l}
 \hat b(g)^{l_1-l}\det(J(g,I))^{-l}&l\leq l_1,\\
 \hat c(g)^{l-l_1}\det(J(g,I))^{-l}&l\geq l_1.\end{array}\right.
\end{equation}
This function satisfies (\ref{Wsharpwelldeflemmaeq2}) with
\begin{equation}\label{l2defeq}
 \chi(\zeta)=\chi_0(\zeta)=\zeta^{l_2},\qquad\text{where}\qquad
 l_2=\left\{\begin{array}{l@{\qquad\text{if }}l}
 l_1-2l&l\leq l_1,\\-l_1&l\geq l_1.\end{array}\right.
\end{equation}
It also satisfies (\ref{Wsharpwelldeflemmaeq4}),
and can therefore be used to define the function $W^\#$ on $G(\R)$ via
\begin{equation}\label{final-arch-W-func-defn}
 W^\#(m_1m_2nk,s)=\delta_P^{s+1/2}(m_1m_2)\chi(m_1)W^{(0)}(m_2)f(k).
\end{equation}
Here, $m_1\in M^{(1)}(\R)$, $m_2\in M^{(2)}(\R)$, $n\in N(\R)$ and $k\in K_\infty$.
It is clear that $W^\#(\,\cdot\,,s)$ satisfies (\ref{Wsharprighttransformeq}), since
$f$ has the corresponding property. By Lemma 2.3.1 of \cite{PS1}, we have
\begin{equation}\label{realWsharpeq}
  W^\#(\eta tuh,s)=\theta(u)^{-1}W^\#(\eta h,s)
\end{equation}
for $t\in T(\R)$, $u\in U(\R)$, $h\in G(\R)$ and
$$
  \eta=\begin{bmatrix}1\\\alpha&1\\&&1&-\bar\alpha\\&&&1\end{bmatrix},
  \qquad \alpha=\frac{b+\sqrt{-D}}{2c},\;D=4ac-b^2.
$$
Note that, if $l=l_1$, then $W^\#$ coincides with the archimedean section
used in \cite{Fu} and \cite{PS1}.
\subsection{The local archimedean integral}
Let $B$ and $W^\#$ be as defined in Sect.\ \ref{realbesselfnsec}
and \ref{realW sharp function}. By (\ref{realBproperty1eq}) and
(\ref{realWsharpeq}), it makes sense to consider the integral
\begin{equation}\label{realZinftyeq2}
 Z(s,W^\#,B)=\int\limits_{R(\R)\backslash H(\R)}W^\#(\eta h, s) B(h)dh.
\end{equation}
Our goal in the following is to evaluate this integral. The
function $W^\#(\eta h, s) B(h)$ is right invariant under $K^H_\infty$. Using this
fact and the disjoint double coset decomposition (\ref{HRBesseldecompeq}), we obtain
\begin{align}\label{archintegral1eq}
 Z(s,W^\#,B)&= \pi\int\limits_{\R^{\times}}\int\limits_1^{\infty}
  W^\#\Big(\eta \begin{bmatrix}\lambda t_0\mat{\zeta}{}{}{\zeta^{-1}}&\\
  &^tt_0^{-1}\mat{\zeta^{-1}}{}{}{\zeta}\end{bmatrix},s\Big) \nonumber \\
 &\hspace{10ex}B\Big(\begin{bmatrix}\lambda t_0\mat{\zeta}{}{}{\zeta^{-1}}&\\
  &^tt_0^{-1}\mat{\zeta^{-1}}{}{}{\zeta}\,\end{bmatrix}\Big)(\zeta-\zeta^{-3})\lambda^{-4}\,d\zeta\,d\lambda;
\end{align}
see (4.6) of \cite{Fu} for the relevant integration formulas. The
above calculations are valid for any choice of $a,b,c$ as long as
$S=\mat{a}{b/2}{b/2}{c}$ is positive definite. To compute
(\ref{archintegral1eq}), we will fix $D = 4ac-b^2$ and make
special choices for $a,b,c$. First assume that \underline{$D\equiv
0 \pmod{4}$}. In this case, let $S(-D) := \mat{\frac D4}{0}{0}{1}$.
Then $\eta =\begin{bmatrix}1&&&\\\frac{\sqrt{-D}}{2}&1&&\\&&1&
\frac{\sqrt{-D}}{2}\\ &&&1\end{bmatrix}$, and we can choose $t_0
=\mat{2^{1/2}D^{-1/4}}{}{}{2^{-1/2}D^{1/4}}$. From
(\ref{archBesselformula2eq}) we have
\begin{equation}\label{archbesselformula3eq}
 B\Big(\begin{bmatrix}\lambda t_0\mat{\zeta}{}{}{\zeta^{-1}}&\\
  &^tt_0^{-1}\mat{\zeta^{-1}}{}{}{\zeta}\,\end{bmatrix}\Big) =
  \renewcommand{\arraystretch}{1.3}\left\{\begin{array}{ll}\lambda^le^{-2\pi \lambda
  D^{1/2}\frac{\zeta^2+\zeta^{-2}}2}& \hbox{ if } \lambda > 0,\\
    0& \hbox{ if } \lambda < 0.
 \end{array}\right.
\end{equation}
Next we rewrite the argument of $W^\#$ as an element of $MNK_{\infty}$,
\begin{align*}
 &\eta \begin{bmatrix}\lambda t_0\mat{\zeta}{}{}{\zeta^{-1}}&\\
  &^tt_0^{-1}\mat{\zeta^{-1}}{}{}{\zeta}\end{bmatrix}\\
 &\qquad=\begin{bmatrix}\lambda \mat{D^{-\frac14}
  \Big(\frac{\zeta^2+\zeta^{-2}}{2}\Big)^{-\frac 12}}{}{}
  {D^{\frac 14}\Big(\frac{\zeta^2+\zeta^{-2}}{2}\Big)^{\frac 12}}&\\
  &\mat{D^{\frac 14}\Big(\frac{\zeta^2+\zeta^{-2}}{2}\Big)^{\frac 12}}{}{}{D^{-\frac14}
  \Big(\frac{\zeta^2+\zeta^{-2}}{2}\Big)^{-\frac12}}
 \end{bmatrix} \\
 &\qquad\times\begin{bmatrix}1&-i\zeta^2&&\\0&1&&\\&&1&0\\&&-i\zeta^2&1\end{bmatrix}
 \mat{k_0}{0}{0}{k_0},\qquad\text{where }
  k_0 = (\zeta^2+\zeta^{-2})^{-1/2}\mat{\zeta^{-1}}{i \zeta}{i\zeta}{\zeta^{-1}} \in \SU(2).
\end{align*}
With $f$ as in (\ref{finalfdefeq}), we have
$$
 f(\mat{k_0}{0}{0}{k_0}) = i^{l+l_2}\Big(\frac{\zeta^2+\zeta^{-2}}2\Big)^{-(l+l_2)}.
$$
Therefore, by (\ref{final-arch-W-func-defn}),
\begin{align}\label{whittakerfirstformulaeq}
 &W^\#\Big(\eta \begin{bmatrix}\lambda t_0\mat{\zeta}{}{}{\zeta^{-1}}&\\
  &t_0^{-1}\mat{\zeta^{-1}}{}{}{\zeta}\end{bmatrix},s\Big) \nonumber \\
 &\qquad=i^{l+l_2} \Big(\frac{\zeta^2+\zeta^{-2}}2\Big)^{-(l+l_2)}\Big|\lambda D^{-\frac12}
  \big(\frac{\zeta^2+\zeta^{-2}}{2}\big)^{-1}\Big|^{3(s+\frac12)}
  \omega_\tau(\lambda)^{-1}
  W^{(0)}(\mat{\lambda D^{\frac12}\big(\frac{\zeta^2+\zeta^{-2}}{2}\big)}{0}{0}{1}).
\end{align}
Let $q\in\C$ be such that $\omega_\tau(y)=y^q$ for $y>0$. It
follows from (\ref{W1Whittakerformulaeq}),
(\ref{archbesselformula3eq}) and (\ref{whittakerfirstformulaeq}) that
\begin{align}\label{realI1eq}
 Z(s,W^\#,B)&=i^{l+l_2}a^+\pi D^{-\frac{3s}2-\frac 34+\frac q4}(4\pi)^{\frac q2}
   \int\limits_0^{\infty}\int\limits_1^{\infty}
   \lambda^{3s+\frac 32+l-\frac q2}
   \Big(\frac{\zeta^2+\zeta^{-2}}2\Big)^{-3s-\frac 32+\frac q2-l-l_2}\nonumber\\
  &\hspace{20ex}W_{\frac{l_1}2, \frac{ir}2}\big(4 \pi \lambda
   D^{1/2}\frac{\zeta^2+\zeta^{-2}}2\big)e^{-2\pi \lambda
   D^{1/2}\frac{\zeta^2+\zeta^{-2}}2}(\zeta-\zeta^{-3})\lambda^{-4}\,d\zeta\,d\lambda.
\end{align}
The substitution $u = (\zeta^2+\zeta^{-2})/2$ leads to
\begin{align*}
 Z(s,W^\#,B)&=i^{l+l_2}a^+ \pi D^{-\frac{3s}2-\frac 34+\frac q4}(4\pi)^{\frac q2}
   \int\limits_1^{\infty}\int\limits_0^{\infty} \lambda^{3s-\frac 32+l-\frac q2}
   u^{-3s-\frac 32+\frac q2-l-l_2}\\
  &\hspace{20ex}W_{\frac{l_1}2, \frac{ir}2}(4 \pi \lambda
   D^{1/2}u) e^{-2\pi \lambda D^{1/2}u}\,\frac{d\lambda}{\lambda}\,du.
\end{align*}
We will first compute the integral with respect to $\lambda$. For
a fixed $u$ substitute $x = 4 \pi \lambda D^{1/2}u$ to get
$$
 Z(s,W^\#,B)=i^{l+l_2}a^+\pi D^{-3s-\frac l2+\frac q2}(4 \pi)^{-3s+\frac 32 -l+q}
 \int\limits_1^{\infty}u^{-6s-2l-l_2+q}\int\limits_0^{\infty}W_{\frac{l_1}2,
 \frac{ir}2}(x)e^{-\frac x2}x^{3s-\frac 32 + l-\frac q2}\,\frac{dx}x\,du.
$$
Using the integral formula for the Whittaker function from
\cite[p. 316]{MOS}, we get
\begin{align}\label{I1integraleq}
 Z(s,W^\#,B)&=i^{l+l_2}a^+\pi D^{-3s-\frac l2+\frac q2}(4\pi)^{-3s+\frac 32 -l+q}\nonumber\\
  &\hspace{10ex}\frac{\Gamma(3s+l-1+\frac{ir}2-\frac q2)
  \Gamma(3s+l-1-\frac{ir}2-\frac q2)}{\Gamma(3s+l-\frac{l_1}2-\frac 12-\frac q2)}
  \int\limits_1^{\infty}u^{-6s-2l-l_2+q} du \nonumber \\
 &=i^{l+l_2}a^+\pi D^{-3s-\frac l2+\frac q2}\,
  \frac{(4\pi)^{-3s+\frac 32 -l+q}}{6s+2l+l_2-q-1}\,
  \frac{\Gamma(3s+l-1+\frac{ir}2-\frac q2)
  \Gamma(3s+l-1-\frac{ir}2-\frac q2)}{\Gamma(3s+l-\frac{l_1}2-\frac12-\frac q2)}.
\end{align}
Here, for the calculation of the $u$-integral, we have assumed
that ${\rm Re}(6s+2l+l_2-q-1)>0$. In all of this we assumed $D\equiv0$ mod $4$.
If $D\equiv3$ mod $4$, one can proceed as in Sect.\ 4.4 of \cite{PS1}.
We summarize the results.
\begin{theorem}\label{archmaintheorem}
 Let $l$ and $D$ be positive integers such that $D\equiv0,3$ mod $4$.
 Let $S(-D)=\mat{D/4}{}{}{1}$ if $D\equiv0$ mod $4$ and $S(-D)=\mat{(1+D)/4}{1/2}{1/2}{1}$
 if $D\equiv3$ mod $4$. Let $B:\:\GSp_4(\R)\rightarrow\C$ be the function defined
 in (\ref{archBesselformula2eq}), and let $W^\#(\,\cdot\,,s)$ be the function
 defined in (\ref{final-arch-W-func-defn}). Let $l_2\in\Z$ be as in (\ref{l2defeq}).
 Then, for ${\rm Re}(6s+2l+l_2-q-1)>0$,
 the local archimedean integral (\ref{localZseq}) is given by
 \begin{equation}\label{archmaintheoremeq1}
  Z(s,W^\#,B)
  =i^{l+l_2}a^+\pi D^{-3s-\frac l2+\frac q2}\,
  \frac{(4\pi)^{-3s+\frac 32 -l+q}}{6s+2l+l_2-q-1}\,
  \frac{\Gamma(3s+l-1+\frac{ir}2-\frac q2)
  \Gamma(3s+l-1-\frac{ir}2-\frac q2)}{\Gamma(3s+l-\frac{l_1}2-\frac12-\frac q2)}.
 \end{equation}
 Here, $q\in\C$ is related to the central character of $\tau$ via
 $\omega_\tau(y)=y^q$ for $y>0$. The number $r\in\C$ is such
 that (\ref{W1Deltapropertyeq}) holds.
\end{theorem}

Note that, if $l\geq l_1$, so that $l_2=-l_1$, the formula in the theorem simplifies to
\begin{equation}\label{archmaintheoremeq2}
 Z(s,W^\#,B)=i^{l-l_1}\frac{a^+}2\pi D^{-3s-\frac l2+\frac q2}\,(4\pi)^{-3s+\frac 32 -l+q}\,
 \frac{\Gamma(3s+l-1+\frac{ir}2-\frac q2)
 \Gamma(3s+l-1-\frac{ir}2-\frac q2)}{\Gamma(3s+l-\frac{l_1}2 + \frac 12-\frac q2)}.
\end{equation}
In particular, for $l=l_1$, we recover Theorem 4.4.1 of \cite{PS1}. We point out that
in our present approach the number $l_1$ (the $\GL_2$ weight)
can be chosen independently of $l$ (the $\GSp_4$ weight), including the case of
different parity.

\section{An application: special values}\label{special-values-section}
Let $\A$ be the ring of adeles of $\Q$. Let $\pi$ be a cuspidal, automorphic representation of $H(\A)$ associated with a holomorphic Siegel cusp form $\Phi$ of degree $2$.
Our local results are strong enough to obtain an integral representation for the
$\GSp_4\times\GL_2$ $L$-function $L(s,\pi\times\tau)$, where $\tau$ is an \emph{arbitrary}
cuspidal, automorphic representation of $\GL_2(\A)$. In the case that $\tau$ comes
from a holomorphic cusp form of the same weight as $\Phi$, but with arbitrary level
and character, we will use the integral representation to obtain a special
$L$-value result that fits
into the general conjecture of Deligne on special values of $L$-functions.
\subsection{Siegel modular forms and Bessel models}\label{siegel modular form}
We would like to apply the theory outlined above to the case where $\pi$ comes
from a holomorphic Siegel modular form of full level. Following \cite{Fu},
we will impose a condition on the Fourier coefficients of this modular form guaranteeing
the existence of a suitable Bessel model for $\pi$.
Let $\Gamma_2 = \SSp_4(\Z)$. For a positive integer $l$ we denote by
$S_l(\Gamma_2)$ the space of Siegel cusp forms of degree $2$ and
weight $l$ with respect to $\Gamma_2$. Let $\Phi \in S_l(\Gamma_2)$ be a Hecke eigenform.
It has a Fourier expansion
$$
 \Phi(Z) = \sum\limits_{S > 0}a(S,\Phi)e^{2 \pi i \tr(SZ)},
$$
where $S$ runs through all symmetric, semi-integral, positive
definite matrices of size two. We shall make the following
assumption\footnote{The ``Assumption 2'' from \cite{Fu} and \cite{PS1}, namly that
$l$ is a multiple of the number of roots of unity in $\Q(\sqrt{-D})$, is no
longer needed in our current approach.} about the function $\Phi$.
\begin{description}
 \item[Assumption:] $a(S,\Phi) \neq 0$ for some $S =
  \mat{a}{b/2}{b/2}{c}$ such that $b^2-4ac = -D < 0$, where $-D$
  is the discriminant of the imaginary quadratic field $\Q(\sqrt{-D})$.
\end{description}
Strong approximation allows for the definition of an adelic function
$\phi= \phi_{\Phi}$ on $H(\A)$ by
\begin{equation}\label{lift of siegel modular form to group}
 \phi(\gamma h_{\infty} k_0) =
 \mu_2(h_{\infty})^l\det(J(h_{\infty},I))^{-l} \Phi(h_\infty \langle I \rangle),
\end{equation}
where $\gamma \in H(\Q)$, $h_{\infty} \in H^{+}(\R)$, $k_0 \in
\prod\limits_{p < \infty}H(\Z_p)$. Here $I = \mat{i}{}{}{i}$, and
$J(g,Z)=CZ+D$ for $g=\mat{A}{B}{C}{D}$ and $Z$ in the Siegel upper half space.
Note that $\phi$ has archimedean weight $(l,l)$ and is a lowest weight vector
with respect to the action of the Lie algebra. The complex conjugate function $\bar\phi$
has weight $(-l,-l)$ and is a highest weight vector; if $\phi$ lies in a space of
automorphic forms realizing a representation $\pi$, then $\bar\phi$ lies in a space
of automorphic forms realizing the contragredient representation $\tilde\pi$.
Let $\psi = \prod\limits_{p}\psi_{p}$ be the character
of $\Q\backslash \A$ which has conductor $\Z_p$ at every finite prime $p$ and such that
$\psi_{\infty}(x) = e^{-2\pi ix}$ for $x \in \R$. Let
\begin{equation}\label{special four coeff matrix}\renewcommand{\arraystretch}{1.2}
 S(-D) = \left\{\begin{array}{l@{\qquad}l}
    \mat{\frac D4}{0}{0}{1}& \hbox{ if } D \equiv 0 \pmod{4}, \\[3ex]
 \mat{\frac{1+D}{4}}{\frac 12}{\frac 12}{1}& \hbox{ if } D \equiv 3 \pmod{4}.
 \end{array}\right.
\end{equation}
Our quadratic extension is $L=\Q(\sqrt{-D})$. Let $T$ be the subgroup of $\GL_2$
defined in (\ref{TFdefeq}) with $S=S(-D)$. Let $\Lambda$ be an ideal class character of $\Q(\sqrt{-D})$, i.e., a character of
$$
 T(\A)/T(\Q)T(\R) \prod\limits_{p<\infty}(T(\Q_p)\cap \GL_2(\Z_p)),
$$
to be specified further below. Note that if we write $\Lambda=\otimes_v\Lambda_v$
with characters $\Lambda_v$ of $L_v^\times$, then $\Lambda_\infty$ is trivial and
$\Lambda_v$ is unramified for each finite $v$. We define the global
Bessel function of type $(S,\Lambda,\psi)$ associated to $\bar{\phi}$ by
\begin{equation}\label{global Bessel model defn}
B_{\bar{\phi}}(h) = \int\limits_{Z_H(\A)R(\Q)\backslash
R(\A)}(\Lambda \otimes \theta)(r)^{-1}\bar{\phi}(rh)dr,
\end{equation}
where $\theta(\mat{1}{X}{}{1})=\psi(\tr(S(-D)X))$. From \cite[(1-17), (1-19), (1-26)]{Su},
we have, for $h_{\infty} \in H^+(\R)$,
\begin{equation}\label{Bessel model arch formula}
 B_{\bar{\phi}}(h_{\infty}) = \mu_2(h_{\infty})^l\,
 \overline{\det(J(h_{\infty},I))^{-l}}\,e^{-2\pi i\,
 \tr(S(-D)\overline{h_{\infty}\langle I\rangle})}
 \sum\limits_{j=1}^{h(-D)} \Lambda(t_j)^{-1}\overline{a(S_j,\Phi)},
\end{equation}
and $B_{\bar{\phi}}(h_{\infty}) = 0$ for $h_{\infty} \not\in H^+(\R)$.
Here, $h(-D)$ is the class number of $\Q(\sqrt{-D})$, the elements
$t_j$, $j = 1,\ldots,h(-D)$, are
representatives of the ideal classes of $\Q(\sqrt{-D})$, and $S_j$,
$j = 1,\ldots,h(-D)$, are representatives of the $\SL_2(\Z)$
equivalence classes of primitive semi-integral positive definite
matrices of discriminant $-D$ corresponding to $t_j$. Thus, by the
Assumption above, there exists a $\Lambda$ such that
$B_{\bar{\phi}}(I_4) \neq 0$. We fix such a $\Lambda$. Let $\mathcal{B}$ be the
space of Bessel functions generated by $B_{\bar\phi}$ under right translation.
Let $\mathcal{B}=\mathcal{B}_1\oplus\ldots\oplus\mathcal{B}_j$ be a decomposition
into irreducible components. Each $\mathcal{B}_i$ has a holomorphic discrete
series representation with scalar minimal $K$-type $(l,l)$ as its archimedean
component, and a spherical representation determined by the Hecke $p$-eigenvalues
as its component at a finite prime $p$. Therefore, all the $\mathcal{B}_i$ are
isomorphic. It follows from the uniqueness of local Bessel models that all spaces
$\mathcal{B}_i$ are identical, i.e., $\mathcal{B}$ is irreducible\footnote{Since
multiplicity one for $\GSp(4)$ is still an issue, we are being
careful here and avoid assuming that $\phi$ itself generates an irreducible
representation.}. If $\pi=\otimes\pi_p$ is any of the irreducible components of the
automorphic representation generated by $\phi$, then the representation of $H(\A)$ on
$\mathcal{B}$ is $\tilde\pi=\otimes\tilde\pi_p$. Since the vector $B_{\bar\phi}$ is
$H(\Z_p)$ invariant for each finite $p$, and a highest weight vector at
the archimedean place, it follows that $B_{\bar\phi}$ is a pure tensor. More precisely,
\begin{equation}\label{Bfactorizationeq}
 B_{\bar\phi}(g)=\overline{a(\Lambda)}
 \prod_{p\leq\infty}B_p(g_p),\qquad g=(g_p)_p\in H(\A),
\end{equation}
where $B_\infty$ is the function given in (\ref{archBesselformula2eq}),
where $B_p$, $p<\infty$, is the spherical vector in the local Bessel model
$\mathcal{B}_{\Lambda_p,\theta_p,\psi_p}(\tilde\pi_p)$, normalized by $B_p(1)=1$
(see Sect.\ \ref{sphericalbesselfunctionsec}),
and where $a(\Lambda) =\sum\limits_{j=1}^{h(-D)} \Lambda(t_j)a(S_j,\Phi)$.
\subsection{Elliptic modular forms as adelic functions}
Let $S_l(N,\chi')$ be the space of  holomorphic cusp forms on the complex upper half plane $\SH_1$ of weight $l$ with respect to $\Gamma_0(N)$ and nebentypus $\chi'$. Here $N = \prod_{p} p^{n_p}$ is any positive integer and $\chi'$ is a Dirichlet character modulo $N$. Then $\Psi \in S_l(N, \chi')$ satisfies
\begin{equation}\label{modular-form-trans-property}
 \Psi\big(\frac{az+b}{cz+d}\big) = \chi'(d) (cz+d)^l \Psi(z) \qquad
  \mbox{ for } z \in \SH_1 \mbox{ and } \mat{a}{b}{c}{d} \in \Gamma_0(N)
\end{equation}
and has a Fourier expansion
\begin{equation}\label{Psifouriereq}
 \Psi(z) = \sum\limits_{n=1}^\infty b_n e^{2 \pi i n z}.
\end{equation}
We will assume that $\Psi$ is primitive, which means that $\Psi$ is a newform,
a Hecke eigenform, and is normalized so that $b_1 = 1$. We will now define a function
$f_\Psi$ on $\GL_2(\A)$ associated to $\Psi$. For this, let
$\omega = \otimes \omega_p$ be the character of $\A^\times / \Q^\times$ defined as the
composition
$$
 \A^\times=\Q^\times\times\R^\times_+\times\Big(\prod_{p<\infty}\Z_p^\times\Big)
 \longrightarrow\prod_{p|N}\Z_p^\times\longrightarrow\prod_{p|N}(\Z_p/p^{n_p}\Z_p)^\times
 \cong(\Z/N\Z)^\times\stackrel{\chi'}{\longrightarrow}\C^\times.
$$
By definition, for primes $p \nmid N$, the local character $\omega_p$ of
$\Q_p^\times$ is unramified and satisfies $\chi'(p) = \omega_p(p)^{-1}$.
Furthermore, $\omega_\infty$ is trivial on $\R^\times_+$ and $\omega_\infty(-1) = (-1)^l$. For primes $p | N$, $\omega_p$ is trivial on $1+p^{n_p}\Z_p$.
Furthermore, for any positive integer $a$ coprime to $N$,
\begin{equation}\label{chi-value-eqn}
 \chi'(a)=\prod\limits_{p | N} \omega_p(a).
\end{equation}
Let $K^{(0)}(N):=\prod\limits_{p|N}K^{(0)}(p^{n_p}\Z_p)\prod\limits_{p\nmid N}\GL_2(\Z_p)$
with the local congruence subgroups
$K^{(0)}(p^n\Z_p)=\GL_2(\Z_p)\cap\mat{1+p^n\Z_p}{\Z_p}{p^n\Z_p}{\Z_p}$
as in (\ref{K'defeq}). Let
$K_0(N):=\prod\limits_{p|N}K_0(p^{n_p}\Z_p)\prod\limits_{p\nmid N}\GL_2(\Z_p)$, where
$K_0(p^n\Z_p)=\GL_2(\Z_p)\cap\mat{\Z_p}{\Z_p}{p^n\Z_p}{\Z_p}$. Evidently,
$K^{(0)}(N)\subset K_0(N)$. Let $\lambda$ be the character of $K_0(N)$ given by
\begin{equation}\label{character-on-K0(N)-def}
 \lambda(\mat{a}{b}{c}{d}) := \prod\limits_{p | N} \omega_p(a_p).
\end{equation}
With these notations, we now define the adelic function $f_\Psi$ by
\begin{equation}\label{fPsidefeq}
 f_{\Psi}(\gamma mk)=\lambda(k)\frac{\det(m)^{l/2}}{(\gamma i
 + \delta)^l}\Psi\Big(\frac{\alpha i + \beta}{\gamma i +\delta}\Big),
\end{equation}
where $\gamma\in\GL_2(\Q)$, $m=\mat{\alpha}{\beta}{\gamma}{\delta} \in \GL_2^{+}(\R)$
and $k\in K_0(N)$. Using (\ref{modular-form-trans-property}),
(\ref{chi-value-eqn}) and (\ref{character-on-K0(N)-def}), it is easy to check
that $f_\Psi$ is well-defined.
Let $V_\Psi$ be the space of right translates of $f_\Psi$, on which the group $\GL_2(\A)$ acts by right translation to give an irreducible, cuspidal, automorphic representation $\tau=\tau_\Psi$. Note that the central character $\omega_\tau$ of $\tau$ is given by the character $\omega$ associated to $\chi'$.

\begin{lemma}\label{identifylocalvectorslemma}
 Let $\tau=\otimes\tau_p$ be the decomposition of $\tau$ into a restricted tensor
 product of local representations. Let the global character $\psi$
 be as in Sect.\ \ref{siegel modular form}. Consider the function
 $$
  W^{(0)}(g)=\int\limits_{\Q\backslash\A}\psi(x)f_\Psi(\mat{1}{x}{}{1}g)\,dx,
 $$
 which is a vector in the global $\psi^{-1}$ Whittaker model
 $\mathcal{W}(\tau,\psi^{-1})=\otimes\mathcal{W}(\tau_p,\psi_p^{-1})$ corresponding
 to the automorphic form $f_\Psi$. This function is a pure tensor of local Whittaker
 functions,
 \begin{equation}\label{identifylocalvectorslemmaeq1}
  W^{(0)}(g)=\prod_{p\leq\infty}W^{(0)}_p(g_p),\qquad g=(g_p)\in\GL_2(\A).
 \end{equation}
 For each finite $p$, the function $W^{(0)}_p$ is the local newform in
 $\mathcal{W}(\tau_p,\psi_p^{-1})$ described in Lemma \ref{GL2-newform-Kirillov}.
 If the normalization is such that $W^{(0)}_p(1)=1$ for all finite $p$, then
 \begin{equation}\label{identifylocalvectorslemmaeq2}
  W^{(0)}_\infty(\mat{t}{0}{0}{1})=\left\{\begin{array}{l@{\qquad\text{if }}l}
  e^{-2\pi t}t^{l/2}&t>0,\\0&t<0.\end{array}\right.
 \end{equation}
 The constants appearing in (\ref{W1Whittakerformulaeq}) are
 $a^+=(4\pi)^{-l/2}b_1=(4\pi)^{-l/2}$ and $a^-=0$. Furthermore,
 \begin{equation}\label{identifylocalvectorslemmaeq3}
  W^{(0)}_\infty(g\mat{\cos(\theta)}{\sin(\theta)}{-\sin(\theta)}{\cos(\theta)})
  =e^{il\theta}W^{(0)}_\infty(g)\qquad\text{for all }g\in\GL_2(\R),\:\theta\in\R.
 \end{equation}
\end{lemma}
{\bf Proof.} By definition, the function $f_\Psi$, and hence the function $W^{(0)}$,
is right invariant under $K^{(0)}(p^{n_p}\Z_p)$ for each finite $p$. Our requirement
that $\Psi$ is a newform implies that $n_p$ is the conductor of the local representation
$\tau_p$. This implies that $W^{(0)}$ factors as in (\ref{identifylocalvectorslemmaeq1}),
and that $W^{(0)}_p$ is the local newform, for each finite $p$.
Formula (\ref{identifylocalvectorslemmaeq2}) follows from a standard calculation.
The equality $a^+=(4\pi)^{-l/2}$ follows from the fact that for the classical Whittaker
function in (\ref{W1Whittakerformulaeq}) we have $ir=l-1$ (for the discrete series
representation of lowest weight $l$ under consideration) and the explicit formula
\begin{equation}\label{identifylocalvectorslemmaeq4}
 W_{\frac l2,\frac{l-1}2}(x)=e^{-x/2}x^{l/2}\qquad\text{for all }x>0.
\end{equation}
The property (\ref{identifylocalvectorslemmaeq3}) follows from (\ref{fPsidefeq}).\qed
\subsection{Choosing the global characters}
In the following we will make a choice for the section $f_\Lambda(g,s)\in I_\C(s,\chi,\chi_0,\tau)$
appearing in the global integrals (\ref{globalintegraleq}) and (\ref{basicidentityeq}).
We will choose $f_\Lambda$ as a pure tensor $\otimes f_p$ via the middle isomorphism in
(\ref{locglobindrepdiagrameq}). Each $f_p$ will be chosen to be the local section
corresponding to the function $W^\#$ from the local integral representations obtained in
Theorem \ref{ram-central-char-main-thm-local} (non-archimedean case)
and Theorem \ref{archmaintheorem} (archimedean case).
We have to make sure, however, that the local data entering these theorems,
in particular the characters $\chi$, $\chi_0$ and $\Lambda$, fit into a global situation.

\begin{lemma}\label{globalchi0lemma}
 Let $L$ be an imaginary quadratic field extension of $\Q$.
 Let $\omega=\otimes\omega_p$ be a character of $\Q^\times\backslash\A^\times$.
 Let $l_2$ be an integer such that $(-1)^{l_2}=\omega_\infty(-1)$. Then there exists
 a character $\chi_0=\otimes\chi_{0,v}$ of $L^\times\backslash\A_L^\times$ such that
 \begin{enumerate}
  \item the restriction of $\chi_0$ to $\A^\times$ coincides with $\omega$, and
  \item $\chi_{0,\infty}(\zeta)=\zeta^{l_2}$ for all $\zeta\in S^1$.
 \end{enumerate}
\end{lemma}
{\bf Proof.} Since $\omega$ is trivial on $L^\times\cap\A^\times=\Q^\times$, we
can extend $\omega$ to a character of $L^\times\A^\times$ in such a way that
$\omega\big|_{L^\times}=1$. Since $S^1\cap(L^\times\A^\times)=\{\pm1\}$, we
can further extend $\omega$ to a character of $S^1L^\times\A^\times$ in such
a way that $\omega(\zeta)=\zeta^{l_2}$ for all $\zeta\in S^1$.
For each finite place $v$ of $L$ we will choose a compact subgroup $U_v$ of $\OF_{L,v}^\times$
such that $\omega$ can be extended to $S^1L^\times\A^\times\big(\prod_{v<\infty}U_v\big)$, with $\omega$ trivial on $\prod_{v<\infty}U_v$ and $U_v=\OF_{L,v}^\times$ for almost all $v$. Hence, the $U_v$ should be chosen such that $\omega$ is trivial on $\big(\prod_{v<\infty}U_v\big)\cap S^1L^\times\A^\times$. We consider the intersection
\begin{equation}\label{globalchi0lemmaeq1}
 \big(\prod_{v<\infty}U_v\big)\cap S^1L^\times\A^\times
 =\big(\prod_{v<\infty}U_v\big)\cap\C^\times L^\times\big(\prod_{p<\infty}\Z_p^\times\big).
\end{equation}
Let $z\alpha x$ be an element of this intersection, where $z\in\C^\times$,
$\alpha\in L^\times$ and $x\in\prod_{p<\infty}\Z_p^\times$. We have
$\alpha\in L^\times\cap\prod_{v<\infty}\OF_{L,v}^\times=\OF_L^\times$, which is a finite set, say $\{\alpha_1,\ldots,\alpha_m\}$.
For $i$ such that $\alpha_i\notin\Q$, choose a prime $p$ such that
$\alpha_i\notin\Z_p^\times$. Then choose a place $v$ lying above $p$, and choose $U_v$ so
small that $\alpha_i\notin U_v\Z_p^\times$. Then the intersection (\ref{globalchi0lemmaeq1})
equals
\begin{equation}\label{globalchi0lemmaeq2}
 \big(\prod_{v<\infty}U_v\big)\cap\C^\times\Q^\times\big(\prod_{p<\infty}\Z_p^\times\big).
\end{equation}
We can choose $U_v$ even smaller, so that $\omega$ is trivial
on this intersection. We can therefore extend $\omega$ to a character of
\begin{equation}\label{globalchi0lemmaeq3}
 S^1L^\times\A^\times\big(\prod_{v<\infty}U_v\big)
 =\C^\times L^\times\big(\prod_{v<\infty}U_v\big)\big(\prod_{p<\infty}\Z_p^\times\big).
\end{equation}
in such a way that $\omega$ is trivial on $\prod_{v<\infty}U_v$.
The group (\ref{globalchi0lemmaeq3}) is of
finite index in $\C^\times L^\times\big(\prod_{v<\infty}\OF_{L,v}^\times\big)$, and therefore
of finite index in $\A_L^\times$ (using the finiteness of the class number).
By Pontrjagin duality, we can now extend $\omega$ to a character $\chi_0$ of
$\A_L^\times$ with the desired properties.\qed

We apply this lemma with $\omega=\omega_\tau$, the central character of the
representation $\tau$ of $\GL_2(\A)$ generated by the cusp form $\Psi$,
and with $l_2=-l$.
We let $\chi_0=\otimes\chi_{0,v}$ be a character of $L^\times\backslash\A_L^\times$
satisfying properties i) and ii) of the lemma. Let $\chi$ be the character
of $L^\times\backslash\A_L^\times$ defined by
\begin{equation}\label{Lambdachichi0eq}
 \chi(\zeta)=\Lambda(\bar\zeta)^{-1}\chi_0(\bar\zeta)^{-1},\qquad \zeta\in\A_L^\times,
\end{equation}
where $\Lambda$ is as in Sect.\ \ref{siegel modular form}.
Since $\Lambda_\infty$ is trivial, we have
\begin{equation}\label{Lambdachichi0eq2}
 \chi_\infty(\zeta)=\chi_{0,\infty}(\zeta)=\zeta^{l_2}\qquad\text{ for all }\zeta\in S^1.
\end{equation}
\subsection{Defining the global section}
Let the characters $\chi$ and $\chi_0$ be chosen as in the previous section.
Let $f_\Psi$ be the function defined in (\ref{fPsidefeq}), and let $W^{(0)}$ be the
corresponding Whittaker function as in Lemma \ref{identifylocalvectorslemma}.
We extend $W^{(0)}$ to a function on $\GU(1,1;L)(\A)$ via
$$
 W^{(0)}(\zeta g) = \chi_0(\zeta)W^{(0)}(g) \qquad \mbox{ for }
 \zeta \in \A_L^\times,\;g \in \GL_2(\A).
$$
If $W^{(0)}_p$ are the local components of $W^{(0)}$ as in
(\ref{identifylocalvectorslemmaeq1}), then the extended function factors again
as $W^{(0)}(g)=\prod_{p\leq\infty}W^{(0)}_p(g_p)$ with local functions extended to
$\GU(1,1;L_p)(\Q_p)\cong M^{(2)}(\Q_p)$ via
$$
 W^{(0)}_p(\zeta g)=\chi_{0,p}(\zeta)W^{(0)}_p(g)\qquad\text{for }
 \zeta\in L_p^\times,\;g\in\GL_2(\Q_p).
$$
Here, $L_p$ is the quadratic algebra $L\otimes_\Q\Q_p$. Let $s$ be a complex parameter.
For each finite prime $p$, let $W^\#_p(\,\cdot\,,s)$ be the element of
$I_{W_p}(s,\chi_p,\chi_{0,p},\tau_p)$ defined at the beginning of Sect.\ \ref{padiczetasec};
see in particular (\ref{Wsharpformulaeq}). Recall that the support of
$W^\#_p(\,\cdot\,,s)$ is contained in
$M(\Q_p)N(\Q_p)\eta H(\Z_p)\Gamma((p \OF_{L_p})^{n_p})$, where
$\OF_{L_p}$ is the ring of integers in
$L_p$, and where $\Gamma$ denotes a principal congruence subgroup
as in (\ref{princ-cong-defn}). The element $\eta$ is defined in (\ref{basicidentityeq}).
For the archimedean place we define
$W^\#_\infty(\,\cdot\,,s)$ as in Sect.\ \ref{realW sharp function}.
Since we are considering the case $l_1=l$, the function $f$ in (\ref{finalfdefeq})
simplifies to $f(g)=\det(J(g,I))^{-l}$ for all $g\in K_\infty$. Hence, the
formula for $W^\#_\infty(\,\cdot\,,s)$ is
\begin{equation}\label{final-arch-W-func-defn2}
 W^\#_\infty(m_1m_2nk,s)=\delta_P^{s+1/2}(m_1m_2)\chi_\infty(m_1)
 W^{(0)}_\infty(m_2)\det(J(k,I))^{-l},
\end{equation}
where $m_1\in M^{(1)}(\R)$, $m_2\in M^{(2)}(\R)$, $n\in N(\R)$ and $k\in K_\infty$.
The local functions $W^\#_p(\,\cdot\,,s)$ for all places $p$ define a global function
\begin{equation}\label{globalWsharpdefeq}
 W^\#(g,s):=\prod_{p\leq\infty}W^\#_p(g_p,s),\qquad g=(g_p)\in\GU(2,2;L)(\A).
\end{equation}
Hence $W^\#(\,\cdot\,,s)$ is an element of the global induced representation
$I_W(s,\chi,\chi_0,\tau)$.

\begin{lemma}\label{globalWsharppropertieslemma}
 The function $W^\#(\,\cdot\,,s)$ has the following properties.
 \begin{enumerate}
  \item Let $\eta$ be the element of $G(\Q)$ defined in (\ref{basicidentityeq}), and let
   $\eta_N$ be the element of $G(\A)$ whose $p$-component is $\eta$ for $p|N$ and
   $1$ for $p\nmid N$. Then
   $$
    W^\#(g,s) = 0\qquad\text{if }g \not\in 
    M(\A)N(\A)\eta K_{\infty}K^\#_G(N)=M(\A)N(\A)\eta_NK_{\infty}K^\#_G(N).
   $$
  \item $W^\#(\,\cdot\,,s)$ is right invariant under the compact group
   $$
    K_G^\#(N) = \prod\limits_{p | N} H(\Z_p) \Gamma((p \OF_{L_p})^{n_p}) 
    \prod\limits_{p \nmid N} G(\Z_p).
   $$
  \item We have
   \begin{equation}\label{globalWsharppropertieslemmaeq1}
    W^\#(gk,s)=\det(J(k,I))^{-l}W^\#(g,s)\qquad\text{for all }g\in G(\A),\;k\in K_\infty.
   \end{equation}
  \item If $m = m_1m_2$, $m_i \in M^{(i)}(\A)$, $n \in N(\A)$, $k =
   k_0k_{\infty}$, $k_0 \in K_G^\#(N)$, $k_{\infty} \in K_{\infty}$, then
   \begin{equation}\label{globalWsharppropertieslemmaeq2}
    W^\#(mn\eta_Nk,s) = \delta_P^{\frac 12 + s}(m)\chi(m_1)
    \det(J(k_{\infty},I))^{-l}W^{(0)}(m_2).
   \end{equation}
   Recall that $\delta_P(m_1m_2) =|N_{L/\Q}(m_1)\mu_1(m_2)^{-1}|^3$.
 \end{enumerate}
\end{lemma}
{\bf Proof.} i), ii) and iii) are immediate from properties of the local functions.
iv) follows from (\ref{Wsharpformulaeq}) and (\ref{final-arch-W-func-defn2}).\qed

Now let $f_\Lambda(\,\cdot\,,s):\:G(\A)\rightarrow\C$ be the element of $I_\C(s,\chi,\chi_0,\tau)$
corresponding to $W^\#(\,\cdot\,,s)$. By (\ref{fWfreleq}),
\begin{equation}\label{fWfreleq2}
 f_\Lambda(g,s)=\sum_{\lambda\in\Q^\times}W^\#
 \Big(\begin{bmatrix}1&&&\\&\lambda\\&&\lambda\\&&&1\end{bmatrix}g,s\Big),
 \qquad g\in G(\A).
\end{equation}

\begin{lemma}\label{globalfpropertieslemma}
 The function $f_\Lambda(\,\cdot\,,s)$ has the following properties.
 \begin{enumerate}
  \item $f_\Lambda(g,s) = 0$ if $g \not\in M(\A)N(\A)\eta_NK_{\infty}K^\#_G(N)$.
  \item $f_\Lambda(\,\cdot\,,s)$ is right invariant under the compact group
   $K_G^\#(N)$.
  \item We have
   \begin{equation}\label{globalfpropertieslemmaeq1}
    f_\Lambda(gk,s)=\det(J(k,I))^{-l}f_\Lambda(g,s)\qquad\text{for all }g\in G(\A),\;k\in K_\infty.
   \end{equation}
  \item If $m = m_1m_2$, $m_i \in M^{(i)}(\A)$, $n \in N(\A)$, $k =
   k_0k_{\infty}$, $k_0 \in K_G^\#(N)$, $k_{\infty} \in K_{\infty}$, then
   \begin{equation}\label{globalfpropertieslemmaeq2}
    f_\Lambda(mn\eta_Nk,s) = \delta_P^{\frac 12 + s}(m)\chi(m_1)
    \det(J(k_{\infty},I))^{-l}f_\Psi(m_2).
   \end{equation}
 \end{enumerate}
\end{lemma}
{\bf Proof.} i), ii) and iii) follow from the corresponding properties of the
function $W^\#(\,\cdot\,,s)$ given in Lemma \ref{globalWsharppropertieslemma}.

iv) Using (\ref{globalWsharppropertieslemmaeq2}), we calculate
\begin{align*}
 f_\Lambda(mn\eta_Nk,s)&=\sum_{\lambda\in\Q^\times}W^\#\Big(m_1
  \begin{bmatrix}1&&&\\&\lambda\\&&\lambda\\&&&1\end{bmatrix}m_2n\eta_Nk,s\Big)\\
 &=\sum_{\lambda\in\Q^\times}\delta_P^{\frac 12 + s}
  (\begin{bmatrix}1&&&\\&\lambda\\&&\lambda\\&&&1\end{bmatrix}m)\chi(m_1)
    \det(J(k_{\infty},I))^{-l}W^{(0)}(\mat{\lambda}{}{}{1}m_2)\\
 &=\delta_P^{\frac 12 + s}(m)\chi(m_1)
    \det(J(k_{\infty},I))^{-l}\sum_{\lambda\in\Q^\times}W^{(0)}(\mat{\lambda}{}{}{1}m_2)\\
 &=\delta_P^{\frac 12 + s}(m)\chi(m_1)
    \det(J(k_{\infty},I))^{-l}f_\Psi(m_2).
\end{align*}
This concludes the proof.\qed
\subsection{The global integral representation}
Observing (\ref{Zfactorizationeq}), (\ref{localzetaeq})
and (\ref{unramifiedlocalzetaeq}), as well as the local Theorems
\ref{ram-central-char-main-thm-local} and \ref{archmaintheorem}, we now obtain
the following result.
\begin{theorem}\label{globalintreptheorem}
 Let $\Phi \in S_l(\Gamma_2)$ be a Hecke eigenform satisfying the Assumption
 made in Sect.\ \ref{siegel modular form}. Let $\phi$ be the adelic function
 corresponding to $\Phi$, and let $\pi$ be an irreducible component
 of the cuspidal, automorphic representation generated by $\phi$. Let $\tau$ be
 the irreducible, cuspidal, automorphic representation of $\GL_2(\A)$ generated
 by a primitve cusp form $\Psi\in S_l(N,\chi')$, where $N$ is a positive
 integer and $\chi'$ is a Dirichlet character modulo $N$. Let
 the global characters $\chi$, $\chi_0$ and $\Lambda$, as well as the
 global section $f_\Lambda \in I_\C(s,\chi,\chi_0,\tau)$, be chosen as above. Then the
 global integral (\ref{globalintegraleq}) is given by
 \begin{equation}\label{globalintreptheoremeq1}
  Z(s,f_\Lambda,\bar\phi)=\Big(\prod_{p\leq\infty}Y_p(s)\Big)
   \frac{L(3s+\frac12,\pi\times\tilde\tau)}
   {L(6s+1,\omega_\tau^{-1})L(3s+1,\tilde\tau\times\mathcal{AI}(\Lambda))}
 \end{equation}
 with
 \begin{equation}\label{globalintreptheoremeq2}
  Y_\infty(s)=\overline{a(\Lambda)}\pi D^{-3s-\frac l2}\,
   (4\pi)^{-3s+\frac 32 -\frac{3l}2}\,\frac{\Gamma(3s+\frac{3l}2-\frac32)}{6s+l-1}.
 \end{equation}
 Here, $\mathcal{AI}(\Lambda)$ is the automorphic representation of $\GL_2(\A)$
 obtained from $\Lambda$ via automorphic induction. The factor $Y_p(s)$ is one
 for almost all $p$ and depends on $\tau_p$; its precise definition is given
 in Theorem \ref{ram-central-char-main-thm-local}.
 The constant $a(\Lambda)$ is defined at the end of Sect.\ \ref{siegel modular form}.
\end{theorem}
{\bf Proof.} Everything follows from the local Theorems
\ref{ram-central-char-main-thm-local} and \ref{archmaintheorem}, but we have to
observe certain normalizations. For example, the global Bessel function $B_{\bar\phi}$
differs from the product of the local Bessel functions $B_p$ (used in the local
theorems) by a constant; see (\ref{Bfactorizationeq}). This explains the factor
$\overline{a(\Lambda)}$ in (\ref{globalintreptheoremeq2}).
Recall that $\chi_0\big|_{\A^\times}=\omega_\tau$ and
$\Lambda\big|_{\A^\times}=\omega_\pi=1$. Therefore, by
(\ref{Lambdachichi0eq}), $\chi\big|_{\A^\times}=\omega_\tau^{-1}$.
It follows that $\tau\times(\chi\big|_{\A^\times})\cong\tilde\tau$.
Substituting this into Theorem \ref{ram-central-char-main-thm-local}, we obtain
the finite Euler factors in (\ref{globalintreptheoremeq1}). By
Theorem \ref{archmaintheorem}, in its simplified version (\ref{archmaintheoremeq2}),
the archimedean Euler factor is given by
\begin{equation}\label{archmaintheoremeq2b}
 Z(s,W^\#_\infty,B_\infty)=i^{l-l_1}
 \frac{a^+}2\pi D^{-3s-\frac l2+\frac q2}\,(4\pi)^{-3s+\frac 32 -l+q}\,
 \frac{\Gamma(3s+l-1+\frac{ir}2-\frac q2)
 \Gamma(3s+l-1-\frac{ir}2-\frac q2)}{\Gamma(3s+l-\frac{l_1}2 + \frac 12-\frac q2)}.
\end{equation}
Here, $q=0$ since the archimedean central character is trivial on $\R_{>0}$.
The number $l_1$, the $\GL_2$ weight, is equal to $l$;
see (\ref{identifylocalvectorslemmaeq3}).
We have $a^+=(4\pi)^{-l/2}$ by Lemma \ref{identifylocalvectorslemma}. Furthermore,
$ir=\pm(l-1)$ for the discrete series representation in question.
Substituting all of these quantities leads to formula (\ref{globalintreptheoremeq2}).\qed

{\bf Remark:} While, for simplicity's sake, we have formulated the theorem above
only for $\tau$ coming from a holomorphic cusp form of the same weight as $\Phi$,
our local theorems are flexible enough to obtain an integral
representation with $\pi$ as above and $\tau$ an \emph{arbitrary} cuspidal,
automorphic representation of $\GL_2(\A)$. In this more general case we would let
$l_1\in\Z$ be \emph{any} of the weights occuring in $\tau$, and let $l_2$ be the integer
defined in (\ref{l2defeq}). Then the hypotheses of Lemma \ref{globalchi0lemma}
are satisfied, so that we obtain global characters $\chi_0$ and $\chi$.
We would define a function $f:\:K_\infty\rightarrow\C$ as in (\ref{finalfdefeq}) and
$W^\#_\infty:\:G(\R)\rightarrow\C$ as in (\ref{final-arch-W-func-defn}). The
non-archimedean sections $W^\#_p$ would be chosen as above. With these choices,
equation (\ref{globalintreptheoremeq1}) still holds, with
$Y_\infty(s)$ replaced by $\overline{a(\Lambda)}$ times the archimedean
local zeta integral given in (\ref{archmaintheoremeq1}).
\subsection{The classical Eisenstein series on $\GU(2,2)$}
For $Z=\mat{\ast}{\ast}{\ast}{z_{22}} \in \HH_2$, let us denote $z_{22}$ by $Z^\ast$.
Let $\hat{Z}=\frac i2(\,^t\!\bar{Z}-Z)$ for $Z \in \HH_2$. Let ${\rm Im}(z)$ denote the
imaginary part of a complex number $z$, and let $I=\mat{i}{}{}{i} \in\HH_2$.
\begin{lem}\label{alternative f-lambda formula-lemma}
 Let $f_\Lambda$ be the function defined in (\ref{fWfreleq2}).
 For any $g \in G^+(\R)$,
   \begin{equation}\label{alternative f-lambda formula}
    f_\Lambda(g\eta_N,s) = \mu_2(g)^l\det(J(g,I))^{-l}\Big(\frac{\det
    \widehat{g\langle I \rangle}}{{\rm Im}(g \langle I
    \rangle)^\ast}\Big)^{3s+\frac 32 - \frac l2}\Psi((g\langle I\rangle)^\ast).
   \end{equation}
\end{lem}
{\bf Proof.} This follows from a direct calculation as in as in Lemma 5.4.1 of \cite{PS1}.
\qed

As in Sect.\ \ref{furusawa-review}, let
$$
 E(g,s;f_\Lambda) = \sum\limits_{\gamma \in P(\Q) \backslash G(\Q)} f_\Lambda(\gamma g, s).
$$
This series is absolutely convergent for ${\rm Re}(s)>\frac12$.
By Lemma \ref{globalfpropertieslemma} iii), this function satisfies
$E(gk,s;f_\Lambda)=\det(J(k,I))^{-l}E(gk,s;f_\Lambda)$ for all $g\in G(\A)$ and $k\in K_\infty$.
It follows that the function on $G(\A)$ given by
$\mu_2(g)^{-l}\det(J(g,I))^lE(g,s;f_\Lambda)$ is right invariant under $K_\infty$.
Since $K_\infty$ is the stabilizer of $I\in\HH_2$, we can define
a function $\mathcal E$ on $\HH_2$ by the formula
\begin{equation}\label{eis ser on siegel half plane}
 \mathcal E(Z,s):= \mu_2(g)^{-l}\det(J(g,I))^l\,
 E\big(g,\frac s3 + \frac l6 - \frac 12;f_\Lambda\big),
\end{equation}
where $g \in G^+(\R)$ is such that $g\langle I \rangle = Z$. The
series that defines $\mathcal E(Z,s)$ is absolutely
convergent for ${\rm Re}(s) > 3 - l/2$ (see \cite{Kl1}). We have
$l\geq 10$ (see \cite{Kl2}), so that one can set $s=0$ to obtain an Eisenstein
series $\mathcal E(Z,0)$ on $\HH_2$. It follows from Lemma
\ref{alternative f-lambda formula-lemma} that $\mathcal E(Z,0)$ is holomorphic.
This Eisenstein series is a modular form of weight $l$ with respect to
$$
 \Gamma^\#_G(N) := G(\Q) \cap G^+(\R) K_G^\#(N).
$$
Its restriction to $\mathfrak{h}_2$ is a modular form of weight $l$ with respect to $\Gamma^\#_G(N) \cap H(\Q) = \SSp_4(\Z)$. We remark that the Eisenstein series constructed
in \cite{PS1} defines, upon restriction to $\mathfrak{h}_2$, a modular form with respect
to a certain congruence subgroup $\Gamma_H^\#(N)$ of level $N$.
The fact that the Eisenstein series $\mathcal E$ obtained above is a modular form with
respect to the full modular group $\SSp_4(\Z)$ is a direct consequence of the choice
of the non-archimedean sections $W^\#_p$ in Sect.\ \ref{non-arch-section}. Let
$$
 \mathcal E(Z,0) = \sum\limits_{\mathcal S \geq 0} b(\mathcal S, \mathcal
 E)e^{2 \pi i \,\tr(\mathcal SZ)}
$$
be the Fourier expansion of $\mathcal E(Z,0)$,
where $\mathcal{S}$ runs through all hermitian half-integral (i.e., $\mathcal S =
\mat{t_1}{\bar{t}_2}{t_2}{t_3},\: t_1, t_3 \in \Z, \:\sqrt{-D}\,t_2 \in
\OF_{\Q(\sqrt{-D})}$) positive semi-definite matrices of size $2\times2$. By \cite{Ha},
\begin{equation}\label{four coeff algeb}
 b(\mathcal S, \mathcal E) \in \bar{\Q} \qquad\mbox{ for any } \mathcal S.
\end{equation}
Here $\bar{\Q}$ denotes the algebraic closure of $\Q$ in $\C$. The following lemma
shows that the global integral $Z(s,f_\Lambda,\bar\phi)$ defined in (\ref{globalintegraleq})
(with $\bar\phi$ replacing $\phi$) can be expressed as the Petersson inner product
of two classical modular forms.
\begin{lem}\label{integral matching formula}
 We have
 $$
  Z\big(\frac l6-\frac 12,f_\Lambda,\bar\phi\big)=\frac12
  \int\limits_{\SSp_4(\Z)\backslash \SH_2}
  \mathcal E(Z,0) \bar{\Phi}(Z) (\det(Y))^{l-3}\,dX\,dY,
 $$
 where $Z =X + iY$.
\end{lem}
{\bf Proof}. The proof follows exactly as in the proof of Lemma 5.4.2 of \cite{PS1}. It is even simpler in this case since $\mathcal E(Z,0)$ is a modular form with respect to $\SSp_4(\Z)$. \qed
\subsection{The special value}
If $\Phi' \in S_l(\Gamma_2)$ is a Hecke eigenform, then the subfield $\Q(\Phi')$ of $\overline{\Q}$, obtained by adjoining the Hecke eigenvalues of $\Phi'$ to $\Q$, is a totally real number field. For a subring $A$ of $\C$, let $S_l(\Gamma_2, A)$ be the space of modular forms whose Fourier coefficients are contained in $A$. By \cite{Mi}, we can assume that $\Phi \in S_l(\Gamma_2, \Q(\Phi))$. Arguing as in the proof of Lemma 5.4.3 of \cite{PS1}, we get the following result.
\begin{lem}\label{integral-is-alg-lem}
 We have
 \begin{equation}\label{integral-is-alg}
  \frac{Z(\frac l6 - \frac 12, f_\Lambda, \bar\phi)}{\langle \Phi,\Phi \rangle} \in \bar{\Q},
 \end{equation}
 where
 $$
  \langle \Phi,\Phi \rangle = \int\limits_{\SSp_4(\Z)\backslash \SH_2}
  \Phi(Z) \bar{\Phi}(Z) (\det(Y))^{l-3}\,dX\,dY.
 $$
\end{lem}
Let $\langle \Psi,\Psi \rangle_1 = (\SL_2(\Z) : \Gamma_1(N))^{-1}
\int\limits_{\Gamma_1(N)\backslash \SH_1}|\Psi(z)|^2y^{l-2}\,dx\,dy$,
where $\Gamma_1(N) := \{\mat{a}{b}{c}{d} \in \Gamma_0(N) : a, d
\equiv 1 \pmod{N} \}$. We have the following generalization of Theorem 4.8.3 of \cite{Fu}.
\begin{thm}\label{special values thm}
 Let $\Phi$ be a cuspidal Siegel eigenform of weight $l$ with respect to
 $\Gamma_2$ satisfying the assumption from Section \ref{siegel modular form} and
 $\Phi \in S_l(\Gamma^{(2)}, \Q(\Phi))$. Let $\Psi \in S_l(N, \chi')$ be a primitive form,
 with $N = \prod p^{n_p}$ any positive integer and $\chi'$ any Dirichlet
 character modulo $N$.
 Let $\pi_\Phi$ and $\tau_\Psi$ be the irreducible, cuspidal, automorphic representations
 of $\GSp_4(\A)$ and $\GL_2(\A)$ corresponding to $\Phi$ and $\Psi$. Then
 \begin{equation}\label{special-values-eqn}
  \frac{L(\frac l2 - 1, \pi_{\Phi} \times \tilde\tau_{\Psi})}{\pi^{5l-8}
  \langle \Phi,\Phi \rangle \langle \Psi,\Psi \rangle_1} \in \bar{\Q}.
 \end{equation}
\end{thm}

{\bf Proof.} By Theorem \ref{globalintreptheorem}, we have
\begin{equation}\label{special values thmeq1}
 Z(\frac l6 - \frac 12, f_\Lambda, \bar\phi) = C \pi^{4-2l} \frac{L(\frac l2 - 1, \pi_{\Phi}
 \times \tilde\tau_{\Psi})}{L(l-2, \omega_\tau^{-1}) L(\frac{l-1}2, \tau_{\Psi} \times 
 \AI(\Lambda))},
\end{equation}
where
$$
 C = \overline{a(\Lambda)}\,D^{-l+\frac 32}\,2^{-4l+6}\,(2l-5)!\,
 \prod\limits_{p | N} Y_p(\frac l6 - \frac 12).
$$
For $p | N$, the $Y_p$ are defined in Theorem \ref{ram-central-char-main-thm-local}.
We first claim that $C \in \overline{\Q}$. Note that, since $\omega_\tau$ is a character of finite order, we obtain $L_p(l-2, \omega_p^{-1}) \in \overline{\Q}$. It follows from an argument as in the proof of Proposition 3.17 of \cite{Rag} that $L_p((l-1)/2, \tau_p \times \AI(\Lambda_p)) \in \bar{\Q}$. This gives $Y_p(\frac l6 - \frac 12) \in \bar{\Q}$ in all cases, except if $\tau_p = \alpha_p \times \beta_p$, with $\alpha_p$ unramified, $\beta_p$ ramified, $\big(\frac{L_p}{p}\big) = 0$ and $\beta_p \chi_{L_p/\Q_p}$ is unramified. In this case,
\begin{equation}\label{Y-in-special-case}
 Y_p(\frac l6 - \frac 12) = \frac{L_p(l-2, \omega_p^{-1})}{1-\Lambda(\varpi_L)
  (\omega_\pi \beta)^{-1}(p) p^{-l/2+1/2}}.
\end{equation}
We have $\omega_\pi \equiv 1$ and $\Lambda(\varpi_L) = \pm 1$. Using the fact that $\alpha(p) \beta(p) =\omega_p(p)\in \bar\Q$ and
$L_p((l-1)/2, \tau_p \times \AI(\Lambda_p)) \in \bar{\Q}$, it can be deduced
from the third line of (\ref{gl2-gl2-gl1-l-fn}) that
$\beta(p) \in \bar\Q$. It follows that the right hand side of
(\ref{Y-in-special-case}) lies in $\bar\Q$. This proves our claim.
Now it is well-known that
\begin{equation}\label{special values thmeq2}
 \frac{L(l-2, \omega_\tau^{-1})}{\pi^{l-2}} \in \bar\Q
\end{equation}
(see, e.g., \cite{N}, VII.2; the adelic $L$-function $L(s,\omega_\tau^{-1})$ coincides
with the Dirichlet $L$-function $L(s,\chi')$, and we have $(-1)^{l-2}=\chi'(-1)$).
Using \cite{Sh}, by the same
argument as in the proof of Theorem 4.8.3 in \cite{Fu}, we get
\begin{equation}\label{special values thmeq3}
 \frac{L(\frac{l-1}2, \tau_{\Psi} \times \AI(\Lambda))}{\pi^{2l-2}
 \langle \Psi,\Psi \rangle_1} \in \bar{\Q}.
\end{equation}
The assertion now follows by combining (\ref{special values thmeq1}),
(\ref{special values thmeq2}) and (\ref{special values thmeq3}) with
Lemma \ref{integral-is-alg-lem}. \qed

In \cite{Fu}, special value results for full level elliptic modular forms are obtained. In \cite{BH},  holomorphic modular forms for full level, a range of weights and all critical values are considered. In \cite{Sa}, certain squarefree levels for both the Siegel cusp form and elliptic cusp form are considered.

\end{document}